\documentclass[11pt]{amsart}
\usepackage{amscd, amsfonts, amsthm}
\usepackage{graphicx}
\title{Homotopy types of the components of spaces of embeddings of compact polyhedra into 2-manifolds}
\author{Tatsuhiko Yagasaki}
\address{Department of Mathematics, Kyoto Institute of Technology, Matsugasaki, Sakyoku, Kyoto 606, Japan}
\email{yagasaki@ipc.kit.ac.jp}

\newtheorem{theorem}{Theorem}[section]
\newtheorem{proposition}{Proposition}[section] 
\newtheorem{corollary}{Corollary}[section] 
\newtheorem{lemma}{Lemma}[section]

\newtheorem*{Prop 2.3'}{Proposition 2.3$'$}

\theoremstyle{definition}

\def \cal {\mathcal}
\def \phi {\varphi}

\renewcommand{\labelenumi}{(\arabic{enumi})}

\begin{document}
\baselineskip 6 mm

\thispagestyle{empty}

\maketitle
\begin{abstract}
Suppose $M$ is a connected PL 2-manifold and $X$ is a compact connected subpolyhedron of $M$ ($X \neq$ 1pt, a closed 2-manifold). 
Let ${\cal E}(X, M)$ denote the space of topological embeddings of $X$ into $M$ with the compact-open topology 
and let ${\mathcal E}(X, M)_0$ denote the connected component of the inclusion $i_X : X \subset M$ in ${\cal E}(X, M)$. 
In this paper we classify the homotopy type of ${\cal E}(X, M)_0$ in term of the subgroup $G = {\rm Im}[{i_X}_\ast : \pi_1(X) \to \pi_1(M)]$. 
We show that 
if $G$ is not a cyclic group and $M \not\cong {\Bbb T}^2$, ${\Bbb K}^2$ then ${\cal E}(X, M)_0 \simeq \ast$, 
if $G$ is a nontrivial cyclic group and $M \not\cong {\Bbb P}^2$, ${\Bbb T}^2$, ${\Bbb K}^2$ then ${\cal E}(X, M)_0 \simeq {\Bbb S}^1$, and 
when $G = 1$, 
if $X$ is an arc or $M$ is orientable then ${\cal E}(X, M)_0 \simeq ST(M)$ and 
if $X$ is not an arc and $M$ is nonorientable then ${\cal E}(X, M)_0 \simeq ST(\tilde{M})$.
Here ${\Bbb S}^1$ is the circle, ${\Bbb T}^2$ is the torus, ${\Bbb P}^2$ is the projective plane and ${\Bbb K}^2$ is the Klein bottle. 
The symbol $ST(M)$ denotes the tangent unit circle bundle of $M$ with respect to any Riemannian metric of $M$ 
and $\tilde{M}$ denotes the orientation double cover of $M$. 
\end{abstract}

%%%%%%%%%%%%%%%%%%%%%%%%%%%% Section 1 %%%%%%%%%%%%%%%%%%%%%%%%%%%%%%%%%%%%%%%%%%

\section{Introduction}

The homotopy types of the identity components of homeomorphism groups of 2-manifolds have been classified in \cite{Ha, Sc, Ya3}. 
In this article we consider the problem of classifying the homotopy types of embedding spaces into 2-manifolds. 

Suppose $M$ is a connected 2-manifold and $X$ is a compact connected subpolyhedron of $M$ with respect to some triangulation of $M$. 
Let ${\cal E}(X, M)$ denote the space of topological embeddings of $X$ into $M$ with the compact-open topology 
and let ${\cal E}(X, M)_0$ denote the connected component of the inclusion $i_X : X \subset M$ in ${\cal E}(X, M)$. 
The purpose of this article is to describe the homotopy type of ${\cal E}(X, M)_0$ 
in term of the subgroup ${i_X}_\ast \pi_1(X) = {\rm Im}[{i_X}_\ast : \pi_1(X) \to \pi_1(M)]$. 

If $X$ is a point of $M$ then ${\cal E}(X, M) \cong M$, and 
if $X$ is a closed 2-manifold then $X = M$ and ${\cal E}(X, M)_0 = {\cal H}(M)_0$, whose homotopy type is already known \cite{Ha, Sc}. 
Below we assume that {\it $X$ is not a point nor a closed 2-manifold.}

The followings are the main results of this paper:

\begin{theorem} 
Suppose ${i_X}_\ast \pi_1(X)$ is not a cyclic subgroup of $\pi_1(M)$. 
\begin{itemize}
\item[\it (1)] ${\cal E}(X, M)_0 \simeq \ast$ \hspace{3pt} if $M \not\cong {\Bbb T}^2, {\Bbb K}^2$.
\item[\it (2)] ${\cal E}(X, M)_0 \simeq {\Bbb T}^2$ if $M \cong {\Bbb T}^2$.
\item[\it (3)] ${\cal E}(X, M)_0 \simeq {\Bbb S}^1$ if $M \cong {\Bbb K}^2$. 
\end{itemize}
\end{theorem} 

\begin{theorem} 
Suppose ${i_X}_\ast \pi_1(X)$ is a nontrivial cyclic subgroup of $\pi_1(M)$. 
\begin{itemize}
\item[\it (1)] ${\cal E}(X, M)_0 \simeq {\Bbb S}^1$ if $M \not\cong {\Bbb P}^2, {\Bbb T}^2, {\Bbb K}^2$. 
\item[\it (2)] ${\cal E}(X, M)_0 \simeq {\Bbb T}^2$ if $M \cong {\Bbb T}^2$. 
\item[\it (3)] Suppose $M \cong {\Bbb K}^2$. 
\begin{itemize}
\item[\it (i)] ${\cal E}(X, M)_0 \simeq {\Bbb T}^2$ if $X$ is contained in an annulus which does not seprate $M$.
\item[\it (ii)] ${\cal E}(X, M)_0 \simeq {\Bbb S}^1$ if $X$ is not the case (i). 
\end{itemize}
\item[\it (4)] Suppose $M \cong {\Bbb P}^2$. 
\begin{itemize}
\item[\it (i)] ${\cal E}(X, M)_0 \simeq SO(3)/{\Bbb Z}_2$ if $X$ is an o.r. circle in $M$. 
\item[\it (ii)] ${\cal E}(X, M)_0 \simeq SO(3)$ if $X$ is not the case (i). 
\end{itemize}
\end{itemize}
\end{theorem} 

Here ${\Bbb S}^1$ is the circle, ${\Bbb T}^2$ is the torus, ${\Bbb P}^2$ is the projective plane and ${\Bbb K}^2$ is the Klein bottle. 
Finally consider the case where $X$ is null homotopic in $M$. 
We choose a Riemannian manifold structure on $M$ and denote by $S(TM)$ the unit circle bundle of the tangent bundle $TM$. 
Let $\tilde{M}$ denote the orientation double cover of $M$. 

\begin{theorem} Suppose ${i_X}_\ast \pi_1(X) = 1$ (i.e., $X \simeq \ast$ in $M$). 
\begin{itemize}
\item[\it (1)] ${\cal E}(X, M)_0 \simeq S(TM)$ if $X$ is an arc or $M$ is orientable. 
\item[\it (2)] ${\cal E}(X, M)_0 \simeq S(T\tilde{M})$ if $X$ is not an arc and $M$ is nonorientable.
\end{itemize}
\end{theorem} 

Since ${\cal E}(X, M)$ is a topological $\ell^2$-manifold \cite[Theorem 1.2]{Ya2}, 
the topological type of ${\cal E}(X, M)_0$ is determined by the homotopy type of ${\cal E}(X, M)_0$ (Theorems 1.1, 1.2, 1.3). 

Main theorems are deduced through the following considerations: 
Section 2 contains some basic facts on 2-manifolds used in this paper.
In Section 3 it is shown that, except a few cases, ${\cal E}(X, M)_0$ is homotopy equivalent to 
the embedding space ${\cal E}(N, M)_0$ of a regular neighborhood $N$ of $X$ into $M$.
Since $X$ is assumed not to be a closed 2-manifold, it follows that $N$ has a boundary and 
admits a core $Y$ which is a wedge (or a one point union) of circles.  
Theorem 1.1 corresponds to the case where $Y$ includes at least two {\em independent} essential circles.
If $Y$ includes only one {\em independent} essential circle, then we have the case of Theorem 1.2. 
In Sections 4 and 5 we discuss how to eliminate {\em dependent} circles from $Y$ without changing the homotopy type of ${\cal E}(Y, M)_0$.
Based on these observations, Theorems 1.1 and 1.2 can be deduced from the homotopy types of homeomorphism groups of 2-manifolds. 

On the other hand, Theorem 1.3 follows from the direct comparison with the unit circle bundle $ST(M)$, and 
this theorem is regarded as the main result in this article. 
In the proof we need a lemma on canonical extension of embeddings of $X$ into a disk, 
which is deduced from the conformal mapping theorem in the complex function theory. 
These are discussed in Section 6. 

In \cite{Ya1} we stated some partial results on homotopy types of embedding spaces of circles, arcs and disks. 
This article provides with a complete answer on this problem. 
In \cite{Ya1} the proof of the arc case depended on some technical arguments using equivariant homotopy equivalences. 
To avoid them, in this article we make a systematic study on naturality and symmetry property of the canonical extensions of embeddings into a disk. 

%%%%%%%%%%%%%%%%%%%%%%% Section 2 %%%%%%%%%%%%%%%%%%%%%%%%%%%%%%

\section{Preliminaries}

Throughout the paper we follow the next conventions: Spaces are assumed to be separable and metrizable, 
and maps are always continuous. ${\rm Fr}_X A$, ${\rm cl}_X A$ and ${\rm Int}_X A$ denote the frontier, closure and interior of a subset $A$ in $X$. 
On the other hand, $\partial M$ and ${\rm Int}\,M$ denote the boundary and interior of a manifold $M$.
The symbol $\cong$ indicates a homeomorphism and $\simeq$ denotes a homotopy equivalence. 
The term orientation preserving (reversing) is abbreviated as o.p. (o.r.).

First we recall some basic facts on the homeomorphism groups of 2-manifolds. 
Suppose $M$ is a 2-manifold and $X$ is a compact subpolyhedron of $M$ (with respect to some triangulation of $M$). 
Let ${\cal H}_X(M)$ denote the group of homeomorphisms $h$ of $M$ onto itself with $h|_X =id$, equipped with the compact-open topology, 
and let ${\cal H}_X(M)_0$ denote the identity component of ${\cal H}_X(M)$. 

\begin{proposition} (\cite{LM, Ya3}) \\
(1) If $M$ is compact, then ${\cal H}_X(M)$ is an $\ell_2$-manifold. \\
(2) If $M$ is noncompact and connected, then ${\cal H}_X(M)_0$ is an $\ell_2$-manifold. 
\end{proposition}

Here $\ell_2$ is the separable Hilbert space consisting of square sumarable real sequences and 
an $\ell_2$-manifold is a separable merizable space which is locally homeomorphic to $\ell_2$.
An ANR is a retract of an open set of a normed space and it has the homotopy type of a CW-complex.  
Every $\ell_2$-manifold is an ANR and its topological type is determined by its homotopy type. (cf. \cite{vM})

The homotopy type of ${\cal H}_X(M)_0$ is classified as follows:   
We use the following notations: ${\Bbb R}^n$ denotes the Euclidean $n$-space, ${\Bbb S}^n$ the $n$-sphere, 
${\Bbb D}^2$ the 2-disk, ${\Bbb T}^2$ the torus, ${\Bbb M}^2$ the M\"obius band, ${\Bbb P}^2$ the projective plane 
and ${\Bbb K}^2$ denotes the Klein bottle.

\begin{proposition} Suppose $M$ is a connected $2$-manifold and $X$ is a compact subpolyhedron of $M$.
\begin{itemize} 
\item[\it (1)] Suppose $M$ is compact \cite{Ha}, \cite[\S3]{Sc}
\begin{itemize}
\item[\it (i)] ${\cal H}_X(M)_0 \simeq SO(3)$ if $(M, X) \cong ({\Bbb S}^2, \emptyset)$, $({\Bbb P}^2, \emptyset)$.
\item[\it (ii)] ${\cal H}_X(M)_0 \simeq {\Bbb T}^2$ if $(M, X) \cong ({\Bbb T}^2, \emptyset)$. 
\item[\it (iii)] ${\cal H}_X(M)_0 \simeq {\Bbb S}^1$ 
if $(M, X) \cong 
({\Bbb D}^2, \emptyset)$, $({\Bbb D}^2, 0)$, $({\Bbb S}^1 \times [0,1], \emptyset)$, $({\Bbb M}, \emptyset)$, $({\Bbb S}^2, 1pt)$, $({\Bbb S}^2, 2pts)$, \\
\hspace*{45mm} $({\Bbb P}^2, 1pt)$ or $({\Bbb K}^2, \emptyset)$. 
\item[\it (iv)] ${\cal H}_X(M)_0 \simeq \ast$ if $(M, X)$ is not the cases (i), (ii) and (iii). 
\end{itemize}
\item[\it (2)] Suppose $M$ is noncompact \cite{Ya3}
\begin{itemize}
\item[\it (i)] ${\cal H}_X(M)_0 \simeq {\Bbb S}^1$ if $(M, X) \cong ({\Bbb R}^2, \emptyset)$, $({\Bbb R}^2, 1pt)$, 
$({\Bbb S}^1 \times {\Bbb R}^1, \emptyset)$, $({\Bbb S}^1 \times [0, 1), \emptyset)$ or $({\Bbb P}^2 \setminus 1pt, \emptyset)$.
\item[\it (ii)] ${\cal H}_X(M)_0 \simeq \ast$ if $(M, X)$ is not the case (i). 
\end{itemize}
\end{itemize}
\end{proposition}

We also note that ${\cal H}_\partial({\Bbb D}), {\cal H}_\partial({\Bbb M}) \simeq \ast$ \cite[Theorem 3.4]{Ep}. 

The next proposition is an assertion on relative isotopies on 2-manifolds \cite[Theorem 3.1]{Ya3}. 

\begin{proposition} 
Suppose $M$ is a connected $2$-manifold and $N$ is a compact $2$-submanifold of $M$. 
If $(M, N)$ satisfies the following conditions, then ${\cal H}(M)_0 \cap {\mathcal H}_N(M) = {\cal H}_N(M)_0$. 
\begin{itemize}
\item[\it (i)] $N$ has no connected component which is a disk, an annulus or a M\"obius band. 
\item[\it (ii)] $cl(M \setminus N)$ has no connected component which is a disk or a M\"obius band. 
\end{itemize}
\end{proposition}

\begin{Prop 2.3'} (Relative version)
Suppose $M$ is a connected $2$-manifold, $N$ is a compact $2$-submanifold of $M$ and $X$ is a nonempty subset of $N$. 
If $(M, N, X)$ satisfies the following conditions, then ${\cal H}_X(M)_0 \cap {\mathcal H}_N(M) = {\cal H}_N(M)_0$.  
\begin{itemize}
\item[\it (i)]
(a) If $H$ is a disk component of $N$, then $\# (H \cap X) \geq 2$. \\
(b) If $H$ is an annulus or M\"obius band component of $N$, then $H \cap X \neq \emptyset$. 
\item[\it (ii)]
(a) If $L$ is a disk component of $cl(M \setminus N)$, then $\# (L \cap X) \geq 2$. \\ 
(b) If $L$ is a M\"obius band component of $cl(M \setminus N)$, then $L \cap X \neq \emptyset$.
\end{itemize}
\end{Prop 2.3'}

Here $\# X$ denotes the cardinal of the set $X$. 
In Proposition 2.3, the conditions (i) and (ii) imply that $M \not\cong {\Bbb S}^2$, ${\Bbb T}^2$, ${\Bbb P}^2$, ${\Bbb K}^2$. 
Hence the condition (i) of \cite[Theorem 3.1]{Ya3} is redundant. 

Next we recall some fundamental facts on embedding spaces into 2-manifolds.
Suppose $M$ is a 2-manifold and $K \subset X$ are compact subpolyhedra of $M$. 
Let ${\cal E}_K(X, M)$ denote the space of embeddings $f : X \hookrightarrow M$ with $f|_K = id$, equipped with the compact-open topology, 
and let ${\cal E}_K(X, M)_0$ denote the connected component of the inclusion $i_X : X \subset M$ in ${\cal E}_K(X, M)$.  

\begin{proposition} 
${\cal E}_K(X, M)$ is an $\ell_2$-manifold \cite{Ya2}. 
\end{proposition}

In the consideration of homotopy types of the embedding space ${\cal E}(X, M)_0$, 
we may always assume that $\partial M = \emptyset$ according to the next lemma.

\begin{lemma}
Let $\overline{M}$ be a 2-manifold obtained from $M$ by attaching a collar $\partial M \times [0, 1]$ along $\partial M$. 
Then the inclusions ${\cal E}(X, {\rm Int}\,M) \subset {\cal E}(X, M) \subset {\cal E}(X, {\rm Int}\,\overline{M}) \subset {\cal E}(X, \overline{M})$ and 
${\cal E}(X, M)_0 \subset {\cal E}(X, {\rm Int}\,\overline{M})_0 \subset {\cal E}(X, \overline{M})_0$ are homotopy equivalences.
\end{lemma}
 
\begin{proof}
%In fact, ${\cal E}(X, M)_0$ contains a PL-embedding $f : X \to M$ with $f(X) \subset {\rm Int}\,M$ and 
%$f$ induces a homeomorphism ${\cal E}(X, M)_0 \cong {\cal E}(f(X), M)_0$. Hence we may assume that $X \subset {\rm Int}\,M$.
%
Using the collar $\partial M \times [0, 1]$ and a boundary collar of $M$, we can find a non-ambient isotopy $h_t : \overline{M} \to \overline{M}$ 
such that $h_0 = id_{\overline{M}}$, $h_t(\overline{M}) \subset {\rm Int}\,\overline{M}$, $h_t(M) \subset {\rm Int}\,M$ ($0 < t \leq 1$) 
and $h_1(\overline{M}) = M$.
Then the homotopy $\phi_t : {\cal E}(X, \overline{M}) \to {\cal E}(X, \overline{M})$, $\phi_t(f) = h_t f$, satisfies 
the conditions that 
$\phi_0 = id$, 
$\phi_t({\cal E}(X, \overline{M})) \subset {\cal E}(X, {\rm Int}\,\overline{M})$, 
$\phi_t({\cal E}(X, M)) \subset {\cal E}(X, {\rm Int}\,M)$ ($0 < t \leq 1$)  
and $\phi_1({\cal E}(X, \overline{M})) \subset {\cal E}(X, M)$. 
Therefore each inclusion mentioned in Lemma 2.1 is a homotopy equivalence with a homotopy inverse $\phi_1$.
\end{proof}

The homeomorphism group ${\cal H}_K(M)_0$ and the embedding space ${\cal E}_K(X, M)_0$ are joined 
by the restriction map $\pi : {\cal H}_K(M)_0 \to {\cal E}_K(X, M)_0$, $\pi(f) = f|_X$. 
In \cite{Ya2} we have investigated some extension property of embeddings of a compact polyhedron into a 2-manifold, 
based upon the conformal mapping theorem. The result is summarized as follows \cite[Theorem 1.1, Corollary 1.1]{Ya2}:  

\begin{proposition} Suppose $\partial M = \emptyset$. 
Then for every $f \in {\mathcal E}_K(X, M)$ there exist a neighborhood ${\mathcal U}$ of $f$ in ${\mathcal E}_K(X, M)$ and 
a map $\phi : {\mathcal U} \to {\mathcal H}_K(M)_0$ such that $\phi(g)f = g$ for each $g \in {\mathcal U}$ and $\phi(f) = id_M$.
\end{proposition}

\begin{corollary} Suppose $\partial M = \emptyset$. \\
(i) The restriction map $\pi : {\cal H}_K(M)_0 \to {\cal E}_K(X, M)_0$ is a principal bundle with the fiber ${\cal G} \equiv \\ 
{\cal H}_{K}(M)_0 \cap {\mathcal H}_X(M)$, where the group ${\mathcal G}$ acts on ${\cal H}_K(M)_0$ by the right composition. \\
(ii) Suppose $K \subset Y \subset X \subset M$ are compact subpolyhedra of $M$. Then 
the restriction map $p : {\cal E}_K(X, M)_0 \to {\cal E}_K(Y, M)_0$, $p(f) = f|_Y$ is 
a locally trivial bundle with fiber ${\cal F} = {\cal E}_K(X, M)_0 \cap {\cal E}_Y(X, M)$. 
\end{corollary}

\begin{proof}
(ii) Since the restriction map $\pi_Y : {\cal H}_K(M)_0 \to {\cal E}_K(Y, M)_0$ is a principal bundle, 
each $f_0 \in {\cal E}_K(Y, M)_0$ admits an open neighborhood ${\cal U}$ and a section $s : {\cal U} \to {\cal H}_K(M)_0$. 
A fiber preserving homeomorphism $\phi : {\cal U} \times {\cal F} \cong p^{-1}({\cal U})$ is defined by $\phi(f, g) = s(f)g$.
\end{proof}

Proposition 2.3 provides with a sufficient condition for the connectivity of the fiber ${\cal G}$.  
In Section 6.2 we investigate some naturality and symmetry properties of the extension map $\phi$ in Proposition 2.5 in the case where $M$ is a disk.

Finally we list some facts on the fundamental groups of 2-manifolds \cite{Ep}. 
For a group $G$ and $g \in G$, let $\langle g \rangle$ denote the cyclic subgroup of $G$ generated by $g$.

\begin{lemma} 
(i)(a) If a simple closed curve $C$ in $M$ is null-homotopic, then it bounds a disk \cite[Theorem 1.7]{Ep}. \\
\hspace*{4mm} (b) If $C_1$ and $C_2$ are disjoint simple closed curves in $M$ and they are homotopic in $M$, then they bounds an annulus \cite{Ep}. \\
(ii) If $M \not\cong {\Bbb P}$, then $\pi_1(M)$ has no torsion elements \cite[Lemma 4.3]{Ep}. \\
(iii) Suppose $N$ is a connected 2-manifold. If there exists $\alpha$, $\beta \in \pi_1(N)$ such that 
$\alpha \beta \neq \beta \alpha$, $\alpha^2 \beta = \beta \alpha^2$ and $\alpha^2 \neq 1$, then $N \cong {\Bbb K}^2$ \cite[Lemma 2.3]{Ep}. \\
(iv) Suppose $M \not\cong {\Bbb T}^2$, ${\Bbb K}^2$. If $\alpha, \beta \in \pi_1(M)$ and $\alpha \beta = \beta \alpha$, 
then $\alpha, \beta \in \langle \gamma \rangle$ for some $\gamma \in \pi_1(M)$ \cite[Lemma 4.3]{Ep}. \\
(v) Suppose $M \not\cong {\Bbb K}^2$, $C$ is a simple closed curve in $M$ which does not bound a disk or a M\"obius band in $M$, 
$x \in C$ and $\alpha \in \pi_1(M, x)$ is represented by $C$. If $\beta \in \pi_1(M, x)$ and $\beta^k = \alpha^\ell$ for some $k$, 
$\ell \in {\Bbb Z} \setminus \{0 \}$, then $\beta \in \langle \alpha \rangle$ \cite[Lemma 3.1]{Ya3}. 
\end{lemma}

\begin{lemma} 
Suppose $M$ is a connected 2-manifold with $\partial M = \emptyset$ and $X$ is a compact connected subset of $M$.  
If $X \simeq \ast$ in $M$, then $X$ has a closed disk neighborhood $D$ in $M$. 
\end{lemma}

\begin{proof}
Any sufficiently small compact connected 2-manifold neighborhood $N$ of $X$ is null homotopic in $M$. 
Each boundary circle $C_i$ of $N$ is also null homotopic in $M$ and bounds a disk $E_i$.  
If $N \subset E_i$ for some $i$, then we set $D = E_i$. 
Otherwise, $N_1 = N \cup \left( \cup_i E_i \right)$ is a closed 2-manifold, so $N_1 = M$. 
Since $N \simeq \ast$ in $M$, it follows that $\pi_1 \, M = 1$, so $M = {\Bbb S}^2$ and $X$ has a disk neighborhood. 
\end{proof}

%%%%%%%%%%%%%%%%%%%%%%%%%%%%%%%%%%%%  Section 3  %%%%%%%%%%%%%%%%%%%%%%%%%%%%%%%%%%%%

\section{Embedding spaces of regular neighborhoods}

Suppose $M$ is a connected 2-manifold with $\partial M = \emptyset$ and $X$ is a compact connected subpolyhedron of $M$ ($X \neq \emptyset$, 1pt). 
Let $N$ be a regular neighborhood of $X$ in $M$. 
By Corollary 2.1 (ii) we have the locally trivial fiber bundle
\[ {\cal F} \equiv {\cal E}(N, M)_0 \cap {\cal E}_X(N, M) \hookrightarrow {\cal E}(N, M)_0 \stackrel{p}{\longrightarrow} {\cal E}(X, M)_0, \ \ p(f) = f|_X. \]

\begin{proposition}
(1) The fiber ${\cal F} = {\cal E}_X(N, M)_0 \simeq \ast$ and the map $p$ is a homotopy equivalence exactly in the following cases: 
\vskip 1mm
(i) $X$ is not an arc nor a circle, \ (ii) $X$ is an arc and $M$ is orientable, \ (iii) $X$ is an o.p. circle.
\vskip 1mm
(2) In the cases \ (iv) $X$ is an arc and $M$ is nonorientable \ and \ (v) $X$ is an o.r. circle,
\vskip 1mm
\noindent there exists a ${\Bbb Z}_2$-action on ${\cal E}(N, M)_0$ for which the map $p$ factors as 
\[ p : {\cal E}(N, M)_0 \stackrel{\pi}{\longrightarrow} {\cal E}(N, M)_0/{\Bbb Z}_2 \stackrel{q}{\longrightarrow} {\cal E}(X, M)_0, \]
where $\pi$ is a double cover and $q$ is a homotopy equivalence.
\end{proposition}

First we prove the next lemma.

\begin{lemma} 
(1)(a) ${\cal H}_X(M)_0 \cap {\cal H}_N(M) = {\cal H}_N(M)_0 \simeq \ast$. (b) ${\cal E}_X(N, M)_0 \simeq \ast$. \\
(2)(a) ${\cal E}(N, M)_0 \cap {\cal E}_X(N, M) = {\cal E}_X(N, M)_0$ in the cases (i), (ii) and (iii) of Proposition 3.1. \\
\hspace*{4mm} (b) ${\cal E}(N, M)_0 \cap {\cal E}_X(N, M) = {\cal E}_X(N, M) \cong {\cal E}_X(N, M)_0 \times {\Bbb Z}_2$ 
in the cases (iv) and (v) of Proposition 3.1.
\end{lemma}

\begin{proof}
(1) (a) Let $h \in {\cal H}_X(M)_0 \cap {\cal H}_N(M)$ and $h_t \in {\cal H}_X(M)_0$ be an isotopy with $h_0 = id_M$ and $h_1 = h$. 
Since each $h_t$ does not interchange the two sides of any edges of $X$, by cutting $M$ along ${\rm Fr}\,X$, 
we can reduces the situation to the case where $X = \partial M$ (a finite union of circles) and $N$ is a collar of $\partial M$. 
Using a boundary collar we can modify $h_t$ so that $h_t|_N = id_N$. This implies that $h \in {\cal H}_N(M)_0$. 

(b) By Corollary 2.1 (i) and (a) we have the locally trivial bundle
\[ {\cal H}_N(M)_0 \hookrightarrow {\cal H}_X(M)_0 \stackrel{\pi}{\longrightarrow} {\cal E}_X(N, M)_0, \ \ \pi(h) = h|_N \]

Since ${\cal H}_N(M)_0, {\cal H}_X(M)_0 \simeq \ast$ (Proposition 2.2), it follows that ${\cal E}_X(N, M)_0 \simeq \ast$. 

(2)(a) Let $f \in {\cal E}(N, M)_0 \cap {\cal E}_X(N, M)$. Below we show that in each case of (i), (ii) and (iii), 

\begin{itemize}
\item[(\#)] the embedding $f$ does not interchange the two sides of each edge of $X$.
\end{itemize}
Then we can easily construct an isotopy $f_t \in {\cal E}_X(N, M)$ such that $f_0 = id_N$ and $f_1 = f$. 
This implies that $f \in {\cal E}_X(N, M)_0$. 

(i) If $X$ is not an arc nor a circle, then $X$ contains a 2-simplex or a triad (a cone over 3 points). Since $f|_X = id_X$, 
$f$ does not interchange the sides of some edge of $X$. Since $X$ is connected, it follows that $f$ does not interchange the sides of any edge of $X$.

(ii) Suppose $X$ is an arc and $M$ is orientable. Then $N$ is a disk. 
Since $f \in {\cal E}(N, M)_0$, there exists an isotopy $f_t \in {\cal E}(N, M)_0$ 
such that $f_0$ is the inclusion $i_N : N \subset M$ and $f_1 = f$.
Choose any point $x \in X$. Then $f_t$ drags the disk $N$ along the loop $f_t(x)$. 
Since $M$ is orientable, this loop is o.p. and $f_0, f_1 : (N, x) \to (M, x)$ define a same orientation at $x$. 
Since $f|_X = id_X$, the embedding $f$ does not interchange the sides of $X$.

(iii) Suppose $X$ is an o.p. circle. Then $N$ is an annulus. 
Choose a point $x \in X$ and an isotopy $f_t \in {\cal E}(N, M)_0$ with $f_0 = i_N$ and $f_1 = f$. 
We note that 
\begin{itemize}
\item[(\#\#)] $f$ does not interchange the two sides of $X$ iff the loop $f_t(x)$ is o.p.
\end{itemize}
This claim is verified by taking a disk neighborhood $D$ of $x$ in $N$ and applying the same argument as in (2)(ii) to the isotopy $f_t|_D$.
If $M$ is orientable, any loop in $M$ is o.p. and (\#) holds.
Below we assume that $M$ is nonorientable.

($\ast$) Suppose $M \not\cong {\Bbb K}^2$. 

$(\ast)_1$ Suppose $X$ does not bound a disk nor a M\"obius band. 
Let $\alpha, \beta \in \pi_1(M, x)$ denote the classes defined by the circle $X$ with some fixed orientation and the loop $f_t(x)$ respectively.
Since $f|_X = id_X$, the isotopy $f_t|_X$ induces a map $F : X \times {\Bbb S}^1 \to M$ such that 
$F|_{X \times 1} = id_X$ and $F|_{x \times {\Bbb S}^1}$ is the loop $f_t(x)$. 
Since $X \times {\Bbb S}^1$ is a torus and $\alpha, \beta \in {\rm Im}\,F_\ast$, it follows that $\alpha \beta = \beta \alpha$.
Hence $\alpha, \beta \in \langle \gamma \rangle$ for some $\gamma \in \pi_1(M, x)$ (Lemma 2.2 (iv)) 
and $\alpha = \gamma^k$ for some $k \in {\Bbb Z} \setminus \{ 0 \}$ since $X$ does not bound a disk. 
Then $\gamma \in \langle \alpha \rangle$ (Lemma 2.2 (v)) and $\beta \in \langle \alpha \rangle$. 
Since $\alpha$ is o.p., so is $\beta$, and (\#) holds. 

$(\ast)_2$ Suppose $X$ bounds a disk or M\"obius band $E$. 
Using the bundle ${\cal H}(M)_0 \to {\cal E}(N, M)_0$, 
we can find an isotopy $\tilde{f}_t \in {\cal H}(M)_0$ such that $\tilde{f}_0 = id_M$ and $\tilde{f}_t|_N = f_t$.
Since $\tilde{f}_1|_X = id_X$ and $M \not\cong {\Bbb S}^2$, ${\Bbb K}$, we have $\tilde{f}_1(E) = E$.
Hence (\#) holds.

$(\ast\ast)$ Suppose $M \cong {\Bbb K}^2$. 
Let $\pi : {\Bbb T} \to {\Bbb K}$ denote the natural double covering and 
let $m$, $\ell$ denote the meridian and o.p. longitude of ${\Bbb K}$ respectively 
($\pi^{-1}(m) = \{ \pm 1\} \times {\Bbb S}^1$, $\pi^{-1}(\ell) = {\Bbb S}^1 \times \{ \pm 1\} \subset {\Bbb T}^2$). 
Since $X$ is an o.p. circle, it follows that $(M, X) \cong ({\Bbb K}, m)$ if $M \setminus X$ is connected, and that 
$(M, X) \cong ({\Bbb K}, \ell)$ if $M \setminus X$ is not connected. 

%$\pi^{-1}(\ell)$ is the union of two longitude $\ell_{\pm} \equiv {\Bbb S}^1 \times \{ \pm 1\}$ in ${\Bbb T}^2$). 

$(\ast\ast)_1$:  
Suppose $(M, X) = ({\Bbb K}, m)$. 
Consider the two meridians $m_{\pm} = \{ \pm 1\} \times {\Bbb S}^1$ in ${\Bbb T}^2$ with the orientation induced from ${\Bbb S}^1$. 
Let $\overline{x}_{\pm} \in m_{\pm}$ denote the point with $\pi(\overline{x}_{\pm}) = x$. 
Take a unique lift $\tilde{f}_t : m_+ \to {\Bbb T}^2$ of the isotopy $f_t|_m : m \to {\Bbb K}$ ($\pi \tilde{f}_t = f_t \pi$) with $\tilde{f}_0 = id_{m_+}$. 
If the loop $f_t(x)$ is o.r., then its lift $\tilde{f}_t(\overline{x}_+)$ is a path from $\overline{x}_+$ to $\overline{x}_-$.
Since $f_1|_m = id_m$, by the definition of $\pi$ it follows that $\tilde{f}_1 m_+ = (m_-)^{-1}$ as loops. 
Hence $m_- \simeq m_+ \simeq \tilde{f}_1 m_+ = (m_-)^{-1}$ in ${\Bbb T}^2$. But this is impossible. 

$(\ast\ast)_2$: 
Suppose $(M, X) = ({\Bbb K}, \ell)$. 
We regard ${\Bbb S}^1$ as the unit circle in ${\Bbb C}$. 
Consider (a) the coverings $\pi : {\Bbb T}^2 \to {\Bbb K}$ and 
$p : {\Bbb S}^1 \times {\Bbb R}^1 \to {\Bbb S}^1 \times {\Bbb S}^1 = {\Bbb T}^2$, $p(z, t) = (z, e^{it})$, 
(b) the covering involutions 
$r : {\Bbb T}^2 \to {\Bbb T}^2$, $r(z, w) = (-z, -\overline{w})$ ($r^2 = id$, $\pi r = \pi$), and 
$\overline{r} : {\Bbb S}^1 \times {\Bbb R}^1 \to {\Bbb S}^1 \times {\Bbb R}^1$, 
$\overline{r}(z, \pi/2 + t) = (-z, \pi/2 -t)$ ($\overline{r}^2 = id$, $p\overline{r} = rp$), and
(c) the covering transformation $\tau : {\Bbb S}^1 \times {\Bbb R}^1 \to {\Bbb S}^1 \times {\Bbb R}^1$, $\tau(z, t) = (z, t + 2 \pi)$ ($p \tau = p$). 

If the loop $f_t(x)$ is o.r., then 
$f_t|_\ell$ has a unique lift $\tilde{f}_t : \ell \to {\Bbb T}^2$ ($\pi \tilde{f}_t = f_t$) 
such that $\tilde{f}_0(\ell) = {\Bbb S}^1 \times \{ 1 \}$ and $\tilde{f}_1(\ell) = {\Bbb S}^1 \times \{ -1 \}$. 
In turn there is a unique lift $\overline{f}_t : \ell \to {\Bbb S}^1 \times {\Bbb R}^1$ of $\tilde{f}_t$ 
such that $\overline{f}_0(\ell) = {\Bbb S}^1 \times \{ 0 \}$. 
Then $\overline{f}_1(\ell) = {\Bbb S}^1 \times \{ (2 k + 1)\pi \}$ for some $k \in {\Bbb Z}$. 
Consider the embedding $\phi$, $\psi : \ell \times [0, 1] \to {\Bbb S}^1 \times {\Bbb R}^1 \times [0, 1]$  
defined by $\phi(z, t) =  (\overline{f}_t(z), t)$, $\psi(z, t) = (\tau^k \overline{r} \overline{f}_t(z), t)$. 
It follows that 
(i) $\phi(\ell \times [0, 1])$ separates ${\Bbb S}^1 \times {\Bbb R}^1 \times [0, 1]$ into two components, 
(ii) $\tau^k \overline{r} \overline{f}_0(\ell) = {\Bbb S}^1 \times \{ (2 k + 1) \pi \}$, 
$\tau^k \overline{r} \overline{f}_1(\ell) = {\Bbb S}^1 \times \{ 0 \}$, so $\psi(\ell \times [0, 1])$ meets both components, 
hence (iii) $\phi(\ell \times [0, 1]) \cap \psi(\ell \times [0, 1]) \neq \emptyset$.
This means that $\overline{f}_t(x) = \tau^k \overline{r} \overline{f}_t(y)$ for some $x$, $y \in {\Bbb S}^1$ and $t \in [0, 1]$. 
It follows that $\tilde{f}_t(x) = r \tilde{f}_t(y)$, hence $x \neq y$ and $f_t(x) = f_t(y)$. This contradicts that $f_t \in {\cal E}(N, {\Bbb K})$.  
This completes the proof. 
The above argument is essential since there is a homotopy $h_t : \ell \to {\Bbb K}$ such that $h_0 = h_1 = i_\ell$ and the loop $h_t(x)$ is o.r.

(b) The assertion is verified in the proof of Proposition 3.1 (2).
\end{proof}

\begin{proof}[Proof of Proposition 3.1]
(1) The conclusion follows from Lemma 3.1 (1)(b) and (2)(a).

(2) In the cases (iv) and (v) 
there exists a $h \in {\cal H}_X(N)$ such that $h^2 = id_N$ and $h$ interchanges the two sides of $X$ in $N$. 
The group ${\Bbb Z}_2 \cong \{ id_N, h \}$ acts on ${\cal E}(N, M)_0$ by the right composition $f \cdot h = fh$. 
This ${\Bbb Z}_2$-action preserves the fibers of $p$ so that it induces the factorization $p = q \pi$ and also induces an action on ${\cal F}$. 

We have ${\cal F} = {\cal E}_X(N, M) = {\cal E}_X(N, M)_0 \cup ({\cal E}_X(N, M)_0 \cdot h)$ (a disjoint union) $\cong {\cal E}_X(N, M)_0 \times {\Bbb Z}_2$ 
since ${\cal E}_X(N, M) \subset {\cal E}(N, M)_0$ and 
${\cal E}_X(N, M)_0 = \{ f \in {\cal E}_X(N, M) \, : \, f \text{ preserves the two sides of $X$}\}$. 
We also note that the local trivializations of $p$ given in the proof of Corollary 2.1 (ii) preserve the ${\Bbb Z}_2$-actions. 
These observations imply 
that $\pi$ is a double cover and 
that $q$ is a locally trivial fiber bundle with fiber ${\cal F}/{\Bbb Z}_2 \cong {\cal E}_X(N, M)_0 \simeq \ast$ and so $q$ is a homotopy equivalence.
\end{proof}

%%%%%%%%%%%%%%%%%%%%%%%%%%%%%%%%%%%%  Section 4  %%%%%%%%%%%%%%%%%%%%%%%%%%%%%%%%%%%%

\section{Simplification of embedded polyhedra}

Suppose $X = C_1 \cup C_2 \subset M$ is a wedge of two circles with a wedge point $p$ and let $i_X : X \subset M$ denote the inclusion map. 
The notation $C_1 \simeq_p C_2$ in $M$ means that the curve $C_1$ is homotopic to the curve $C_2$ in $M$ relative to $p$ 
under some parametrizations (or orientations) of $C_1$ and $C_2$. 
The notation $(C_1)^2$ denotes a curve which goes around twice along $C_1$ under a parametrization of $C_1$. 

\begin{lemma} 
Suppose $C_1$ and $C_2$ are essential in $M$ and ${i_X}_\ast \pi_1(X, p)$ is a cyclic subgroup of $\pi_1(M, p)$. \\
(1) If both $C_1$ and $C_2$ are o.p., 
then $C_1 \simeq_p C_2$ in $M$ and $X$ has a compact neighborhood $M_0$ such that $(M_0, C_1, C_2)$ is as in Figure 4.1 (1). \\
(2) If both $C_1$ and $C_2$ are o.r.,
then $C_1 \simeq_p C_2$ in $M$ and $X$ has a compact neighborhood $M_0$ such that $(M_0, C_1, C_2)$ is as in Figure 4.1 (2). \\
(3) If $C_1$ is o.p. and $C_2$ is o.r., 
then $C_1 \simeq_p (C_2)^2$ in $M$ and $X$ has a compact neighborhood $M_0$ such that $(M_0, C_1, C_2)$ is as in Figure 4.1 (3). 
\end{lemma}

\begin{center}
%\fbox{Figure 4.1 (1)\,(2)\,(3)}
\includegraphics[scale=1.0]{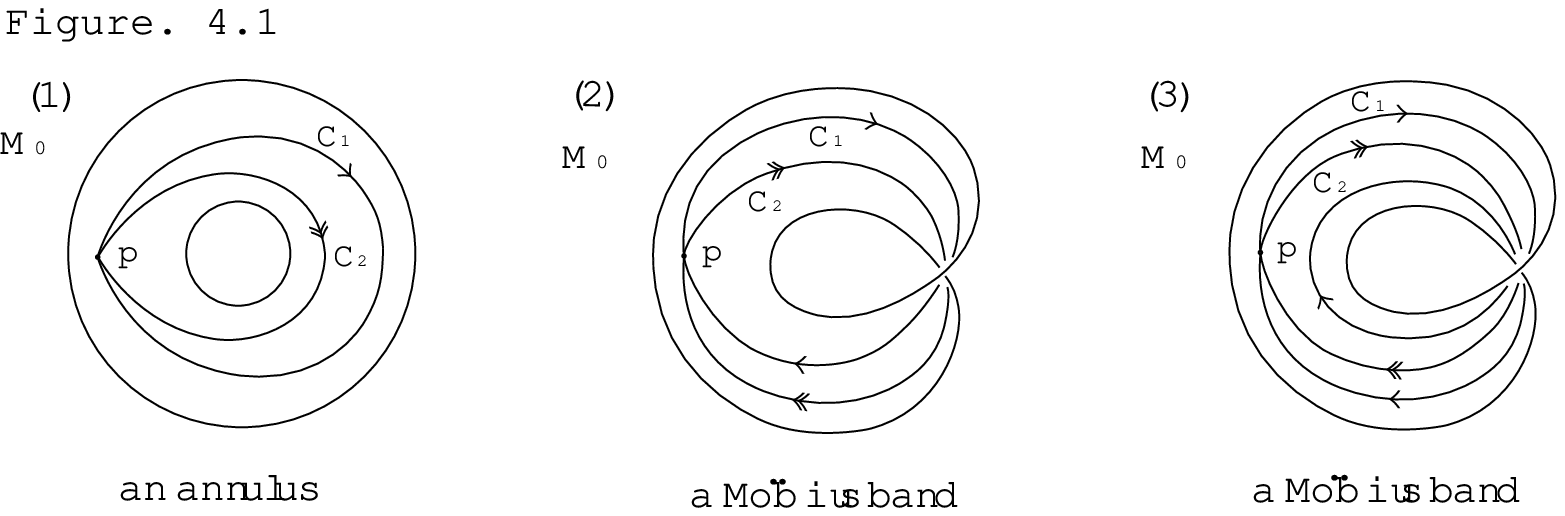} 
\end{center}

\begin{proof} (c.f.\,Figure 4.2) \ 
We choose a small disk neighborhood $N_0$ of $p$ and thin regular neighborhoods $N_1$ and $N_2$ of $C_1$ and $C_2$ 
which are in general position so that $N = N_1 \cup N_2$ is a regular neighborhood of $X$ and $N_1 \cap N_2 = N_0$. 
An outer boundary circle $A$ of $N$ means a boundary circle of $N$ which intersects both $N_1$ and $N_2$.   
The circles $C_1$ and $C_2$ intersect at $p$ either transversely ($\ast$) or tangentially ($\ast\ast$) (Figure 4.2 (0)).  

Let $\alpha$, $\beta \in \pi_1(M, p)$ denote the elements represented by $C_1$ and $C_2$ with some fixed orientations.
Since $\alpha$, $\beta$ are contained in the cyclic subgroup ${i_X}_\ast \pi_1(X, p)$, 
it follows that $\alpha^k = \beta^\ell$ for some $k, \ell \in {\Bbb Z} \setminus \{ 0 \}$.

(1) (cf.\,Figure 4.2 (1)) In this case $N_1$ and $N_2$ are annuli. 
In the case ($\ast$), $N$ is a torus with a hole. 
Let $\widehat{N}$ be the torus obtained from $N$ by attaching a disk along the boundary circle $\partial N$, 
and let $r : M \to \widehat{N}$ be a map with $r|_N = id_N$. 
Then $r_\ast \alpha$ and $r_\ast \beta$ are two generators of $\pi_1(\hat{N}, p) \cong {\Bbb Z} \oplus {\Bbb Z}$, 
while $r_\ast \alpha$ and $r_\ast \beta$ are contained in the cyclic subgroup $(r i_X)_\ast \pi_1(X, p)$. This is a contradiction. 
Therefore, ($\ast\ast$) holds and $N$ is a sphere with three holes (Figure 4.2 (1)(a)).  

We note the following facts: \\
(i) $M \not\cong {\Bbb P}^2$.
In fact, if $M \not\cong {\Bbb P}^2$, 
then $\alpha \neq 0$ is the generator of $\pi_1(M, p) = {\Bbb Z}_2$ and $C_1$ is o.r., a contradiction. \\
(ii) If $X$ is contained in an annulus or a M\"obius band $E$ in $M$, 
then $A$ bounds a disk $D$ in $E$ and $M_0 = N \cup D$ satisfies the required condition. \\
(iii) Consider the case where $M \cong {\Bbb K}$. 
If $C_1$ does not separate $M$, then $X$ is contained in an annulus. 
If $C_1$ separates $M$, then $C_1$ is a common boundary circle of two M\"obius bands, one of which contains $X$. 
In either case the required conclusion follows from (ii). 

Below we assume that  
\vskip 1mm
\hspace*{1mm} $(\#)$ \hspace{2mm} $M \not\cong {\Bbb K}$ \ and \ $X$ is not contained in any annulus and M\"obius band.
\vskip 1mm
\noindent and derive a contradiction.

First we show that $C_1 \simeq_p C_2$ in $M$ (i.e., $\beta = \alpha^{\pm}$). 
Since $C_1$ and $C_2$ are essential, they do not bound a disk. 
Suppose $C_1$ bounds a M\"obius band $L$. 
By the assumption we have $C_2 \subset cl(M \setminus L)$. 
Let $\overline{M}$ denote the 2-manifold obtained from $M$ by replacing $L$ by a disk $E$ and 
take a map $r : M \to \overline{M}$ such that $r = id$ on $cl(M \setminus L)$ and $r(L) = E$. 
Since $r(C_1) \simeq \ast$ in $\overline{M}$ we have $(r_\ast \beta)^\ell = (r_\ast \alpha)^k = 1$ in $\pi_1(\overline{M})$, 
and since $\overline{M} \not\cong {\Bbb K}$, so $\overline{M} \not\cong {\Bbb P}^2$, we have $r_\ast \beta = 1$ by Lemma 2.2 (ii). 
Hence $r(C_2)$ bounds a disk $F$ in $\overline{M}$ and 
it follows that $r^{-1}(F)$ is a disk bounded by $C_2$ or a M\"obius band bounded by $C_2$ and containing $X$. 
Both cases yield contradictions. Hence $C_1$ does not band a M\"obius band, and similarly $C_2$ does not band a M\"obius band. 
Since $M \not\cong {\Bbb K}$, by Lemma 2.2 (v) we have $\beta \in \langle \alpha \rangle$ and $\alpha \in \langle \beta \rangle$. 
Since $M \not\cong {\Bbb P}$, by Lemma 2.2 (ii) $\pi_1(M)$ has no torsion, so $\beta = \alpha^{\pm}$ as required. 

Now we have the cases (b) and (c) in Figure 4.2 (1), depending on the orientations of $C_1$ and $C_2$. 
In (b), $A \simeq C_1 \ast (C_2)^{-1} \simeq \ast$ and $A$ bounds a disk $D$ in $M$. 
Then $N \cup D$ is an annulus containing $X$, which contradicts $(\#)$. 
In (c), $A_1 \simeq A_2$, so $A_1$ and $A_2$ bounds an annulus B (Lemma 2.2 (i)(b)). 
Let $\overline{M}$ denote the 2-manifold obtained from $M$ by replacing $F = cl(M \setminus (N \cup B))$ by a disk $E$ and 
take a map $r : M \to \overline{M}$ such that $r = id$ on $N \cup B$ and $r(F) \subset E$. 
Then $\overline{M} \cong {\Bbb T}$ or ${\Bbb K}$ and $C_1 \simeq_p (C_2)^{-1}$ in the annulus $N \cup E$.
Since $C_1 \simeq_p C_2$ in $\overline{M}$, it follows that $(C_2)^2 \simeq \ast$. 
Since $\overline{M} \not\cong {\Bbb P}^2$, we have $C_2 \simeq \ast$. This is impossible since $C_2$ is a meridian of $\overline{M}$. 
This completes the proof of (1).

(2) (cf.\,Figure 4.2 (2)) In this case $N_1$ and $N_2$ are M\"obius bands. 
In the case ($\ast\ast$) $N$ is a Klein bottle with a hole. 
Let $\widehat{N}$ be a Klein bottle obtained from $N$ by attaching a disk along the boundary circle 
and let $r : M \to \widehat{N}$ be a map with $r|_N = id_N$. 
It follows that $r_\ast \alpha$ and $r_\ast \beta$ are the center circles of the two M\"obius bands with a common boundary circle. 
However, $r_\ast \alpha$ and $r_\ast \beta$ are contained in the cyclic subgroup $(r i_X)_\ast \pi_1(X, p)$. This is a contradiction. 
Hence the case ($\ast$) holds and $N$ is a M\"obius band with a hole (Figure 4.2 (2)($\ast$)).

We note the following facts: \\
(i) If $M \cong {\Bbb P}^2$, then both $A_1$ and $A_2$ bound disks $D_1$ and $D_2$ respectively and $M_0 = N \cup D_1$ satisfies the required condition. \\
(ii) If $M \cong {\Bbb K}$, then one of $A_1$ and $A_2$ bounds a disk $D$ and another one bounds a M\"obius band. 
The M\"obius band $M_0 = N \cup D$ satisfies the required condition. 

Below we assume that $M \not\cong {\Bbb P}^2$, ${\Bbb K}$. 
Both $C_1$ and $C_2$ do not bound a disk nor a M\"obius band since they are o.r. 
Since $M \not\cong {\Bbb K}$, by Lemma 2.2 (v) 
we have $\beta \in \langle \alpha \rangle$ and $\alpha \in \langle \beta \rangle$. 
Since $M \not\cong {\Bbb P}$, by Lemma 2.2 (ii) $\pi_1(M)$ has no torsion and $\beta = \alpha^{\pm}$, so we have $C_1 \simeq_p C_2$. 
Therefore, we have the situation in Figure 4.2 (2). 
The boundary $\partial N$ consists of two boundary circles $A_1$ and $A_2$,  
one of which is homotopic to $C_1 \ast (C_2)^{-1} \simeq \ast$ and the other is homotopic to $C_1 \ast C_2$. 
The former bounds a disk $D$ and $M_0 = N \cup D$ satisfies the required condition. 

(3) (cf. Figure 4.2 (3)) In this case $N_1$ is an annulus and $N_2$ is a M\"obius band. 
In the case ($\ast$) $N$ is a Klein bottle with a hole. 
Let $\widehat{N}$ be a Klein bottle obtained from $N$ by attaching a disk along the boundary circle, 
and let $r : M \to \widehat{N}$ be a map with $r|_N = id_N$.
Then $r_\ast \alpha$ and $r_\ast \beta$ are two generators of $\pi_1(\widehat{N})$, 
while $r_\ast \alpha$ and $r_\ast \beta$ are contained in the cyclic subgroup $(r i_X)_\ast \pi_1(X, p)$. This is impossible.
Therefore, the case ($\ast\ast$) holds and $N$ is a M\"obius band with a hole (Figure 4.2 (3) ($\ast\ast$)). 

We note the following facts: \\
(i) $M \not\cong {\Bbb P}^2$ as shown in (1)(i) \\
(ii) Suppose $M \cong {\Bbb K}$.
If $C_1$ does not separate $M$, then $X$ is contained in an annulus and $C_2$ is o.p., a contradiction.  
Hence $C_1$ separates $M$ and so $C_1$ is a common boundary circle of two M\"obius bands, one of which contains $X$. 
Then $A$ bounds a disk $D$ in this M\"obius band and $M_0 = N \cup D$ satisfies the required condition.

Below we assume that $M \not\cong {\Bbb K}$. 
If $C_1$ does not bound a M\"obius band, then by Lemma 2.2 (v) we have $\beta \in \langle \alpha \rangle$. 
Since $\alpha$ is o.p., so is $\beta$. This is a contradiction. 
Therefore, $C_1$ bound a M\"obius band $L$. 
Suppose $C_2 \subset cl(M \setminus L)$. 
Let $\overline{M}$ denote the 2-manifold obtained from $M$ by replacing $L$ by a disk $D$  
and let $r : M \to \overline{M}$ denote a map with $r = id$ on $cl(M \setminus L)$ and $r(L) \subset D$. 
Then $(r_\ast \beta)^\ell = r_\ast \alpha^k = 1$. 
Since $\overline{M} \not\cong {\Bbb P}^2$, $\pi_1(\overline{M})$ has no torsion and 
we have $r_\ast \beta = 1$ and $r(C_2)$ is o.p.
By the definition of $r$, this implies that $C_2$ is also o.p., a contradiction. 
Therefore we have $X \subset L$ and the conclusion follows from an easy argumet.
%
%In the case (b) $A_1$ and $A_2$ are freely homotopic in $M$ and they bound an annulus $B$ in $M$. 
%Then it follows that $B$ is the union of $N_2$ and a disk and that the union of $N_2$ and $B$ form a M\"obius band. 
%Then $A_1 \simeq (A_2)^{-2}$, but $A_1 \simeq (A_2)^2$ and this means that $(C_1)^2 \simeq \ast$. Since $C_1 \not\simeq \ast$, 
%it follows that $\pi_1(M)$ has a torsion element and $M$ is a projective plane. Hence $A_1$ bounds a disk in $M$ and hence $C_1 \simeq \ast$, a contradiction. 
%Hence we have the case (a) and it follows that $A \simeq C_1 \ast (C_2)^{-2} \simeq \ast$ in $M$ and that $A$ bounds a disk in $M$. 
%Set $M_0 = N \cup D$.
\end{proof}

%\[ \fbox{Figure 4.2} \]

\vskip 10mm

%\begin{center}
\includegraphics[scale=1.0]{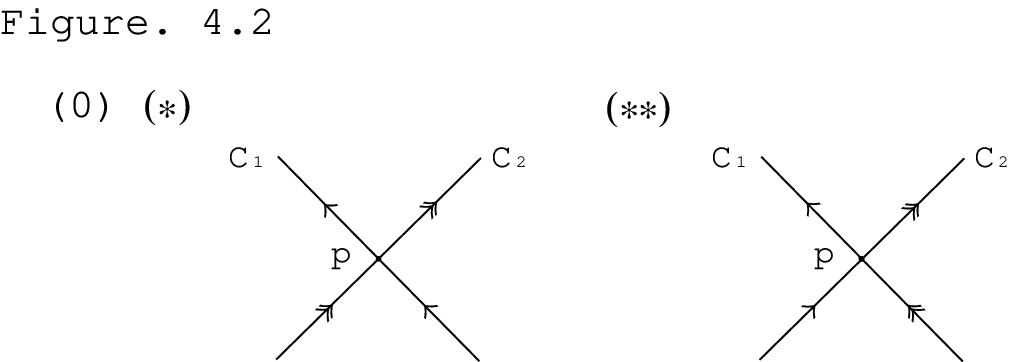} 
%\end{center}

\newpage

%\begin{center}
\hspace*{-5mm}\includegraphics[scale=1.0]{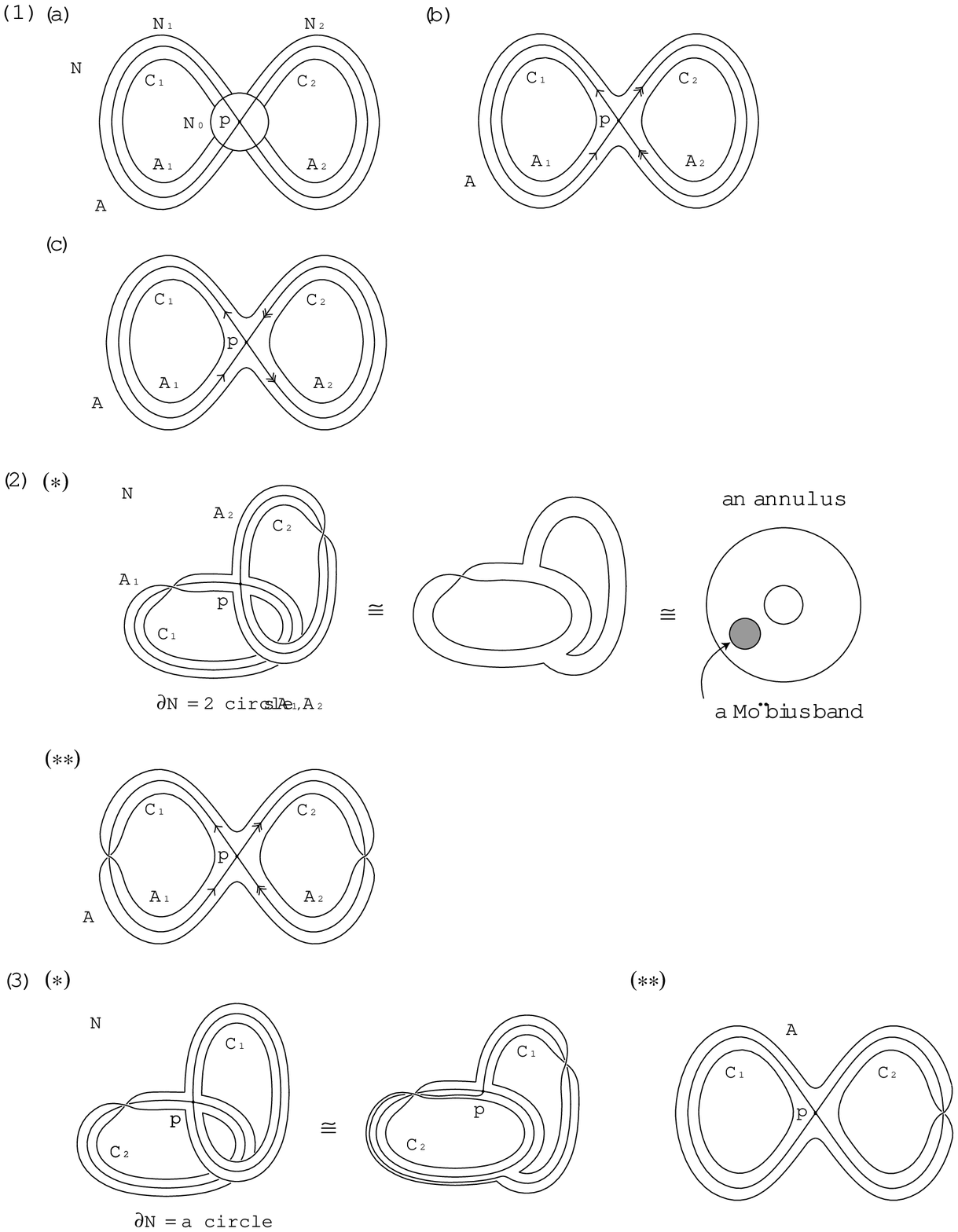} 
%\end{center}

\begin{lemma} Suppose $M$ is a 2-manifold and $X$ is a compact connected subpolyhedron of $M$. 

(1) If $E$ is a disk or a M\"obius band in $M$ and $\partial E \subset X$, 
then the restriction map $p : {\cal E}(X \cup E, M)_0 \to {\cal E}(X, M)_0$ is a homotopy equivalence. 

(2) Suppose $X = Y \cup C_1 \cup C_2$ is a one point union of a compact connected subpolyhedron $Y$ ($\neq$ 1pt) and two essential circles $C_1$ and $C_2$. 
If the pair $(C_1, C_2)$ satisfies one of the conditions listed in Lemma 4.1 (1), (2) and (3), then 
the restriction map $p : {\cal E}(X, M)_0 \to {\cal E}(Y \cup C_2, M)_0$ is a homotopy equivalence. 

(3) Suppose $X = Y \cup C$ is a one point union of a compact connected subpolyhedron $Y$ and a circle $C$. 
If $C \simeq \ast$ in $M$ and $Y$ satisfies one of the conditions (i), (ii) and (iii) in Proposition 3.1, 
then the restriction map $p : {\cal E}(X, M)_0 \to {\cal E}(Y, M)_0$ is a homotopy equivalence. 
\end{lemma}

\begin{proof} (cf.\,Figure 4.3) 

(1) Attaching a collar to $\partial M$, we may assume that $X \cup E \subset {\rm Int}\,M$. Let $N$ denote a regular neighborhood of $X$ in $M$. 
Then $N \cup E$ is a regular neighborhood of $X \cup E$ in $M$. 
Consider the diagram
\[ 
\begin{CD} 
{\cal E}(N \cup E, M)_0 @>p_1>\simeq> {\cal E}(X \cup E, M)_0 \\
@VqVV @VVpV \\
{\cal E}(N, M)_0 @>p_2>\simeq> {\cal E}(X, M)_0.
\end{CD}
\]

By the assumption $X$ is not an arc, and if $X$ is a circle, then $X = \partial E$, which is an o.p. circle.  
Hence by Proposition 3.1 the restriction maps $p_1$ and $p_2$ are homotopy equivalences. 
It suffices to show that the map $q$ is a homotopy equivalence. 
By Corollary 2.1 (ii) the map $q$ forms a locally trivial bundle: 
\[ {\cal E}(N \cup E, M)_0 \cap {\cal E}_N(N \cup E, M) \hookrightarrow {\cal E}(N \cup E, M)_0 \stackrel{q}{\longrightarrow} {\cal E}(N, M)_0 \]
Since (i) $N \cap E = cl(E \setminus \cup_{i=1}^m E_i)$, each $E_i$ is a disk or a M\"obius band in ${\rm Int}\,E$ and they are mutually disjoint and 
(ii) ${\cal H}_\partial(E_i) \simeq \ast$, it follows that 
\[ {\cal E}(N \cup E, M)_0 \cap {\cal E}_N(N \cup E, M) = {\cal E}_N(N \cup E, M)_0 \cong {\cal H}_{N \cap E}(E) \simeq \ast \] 
and hence the map $q$ is a homotopy equivalence as required.

(2) We choose a small disk neighborhood $A$ of the wedge point $x$ of $X$ and 
thin regular neighborhoods $N(Y)$, $N(C_1)$ and $N(C_2)$ of $Y$, $C_1$ and $C_2$.
We may assume that they are in general position and intersects exactly in $A$. 
Thus, for instance, we have that $N(C_1 \cup C_2) = N(C_1) \cup N(C_2)$ is a regular neighborhood of $C_1 \cup C_2$ and 
$N(X) = N(Y) \cup N(C_1) \cup N(C_2)$ is a regular neighborhood of $X$, etc. 
In each case of (1), (2) and (3) in Lemma 4.1, $cl(M_0 \setminus N(C_1 \cup C_2))$ has a unique disk component, which we denote by $D$. 
Note that every component of $cl(D \setminus N(X))$ is a disk. We denote these disk components by $E$, $D_1, \cdots, D_m$ ($m \geq 0$), 
where $E$ is the unique component which meets $N(C_1)$.

Consider the following commutative diagram:
\[ (\ast) \hspace{10mm}
\begin{CD}
{\cal E}(N(X) \cup D, M)_0 @>>> {\cal E}(N(X), M)_0 @>>> {\cal E}(X, M)_0  \\
@V{p_1}VV @. @VV{p}V \\
 {\cal E}(N(Y \cup C_2) \cup \left( \cup_{i=1}^m\, D_i \right), M)_0 @>>> {\cal E}(N(Y \cup C_2), M)_0 @>>> {\cal E}(Y \cup C_2, M)_0.
\end{CD} 
\]
\noindent The horizontal arrows are homotopy equivalences by (1) and Proposition 3.1. Let $F = cl(N(C_1)\setminus A)$. 
\vskip 2mm
It follows that $F$ is a disk, $F \cap E$ is an arc and hence $F \cup E$ is also a disk, and that 
\vskip -3mm
{\small 
\[ 
N(X) \cup D = N(X) \cup \left( \cup_{i=1}^m\, D_i \right) \cup E = \left[ N(Y \cup C_2) \cup F \right] \cup \left[ \left( \cup_{i=1}^m\, D_i \right) \cup E \right] 
= \left[ N(Y \cup C_2) \cup \left( \cup_{i=1}^m\, D_i \right) \right] \cup \left[ F \cup E \right]  
\]
}
\vskip -6mm
\noindent and $F \cup E$ meets the 2-manifold $N(Y \cup C_2) \cup \left( \cup_{i=1}^m\, D_i \right)$ in an arc. 
Therefore,  the map $p_1$ is a homotopy equivalence and the map $p$ is also a homotopy equivalence. 

(3) The proof is essentially same as (2). 
We choose a small disk neighborhood $A$ of the wedge point $x$ of $X$ and thin regular neighborhoods $N(Y)$ and $N(C)$ of $Y$ and $C$. 
The circle $C$ bounds a disk $D_0$ and $D \equiv cl(D_0 \setminus N(C))$ is a disk. Every component of $cl(D \setminus N(X))$ is a disk. 
We denote these components by $E$, $D_1, \cdots, D_m$ ($m \geq 0$), where $E$ is the unique disk component which meets $N(C)$.
Let $F = cl(N(C)\setminus A)$. \
Consider the diagram $(\ast)$ in (2), where $Y \cup C_2$ is replaced by $Y$. 
Then the horizontal arrows and the map $p_1$ are homotopy equivalences by the same reasons, and therefore, the map $p$ is also a homotopy equivalence.
This completes the proof. 
\end{proof}

%\[ \fbox{Figure 4.3 (1), (2), (3)} \]
%\includegraphics[scale=.5]{figure3.eps} 
\begin{center}
\includegraphics[scale=1.0]{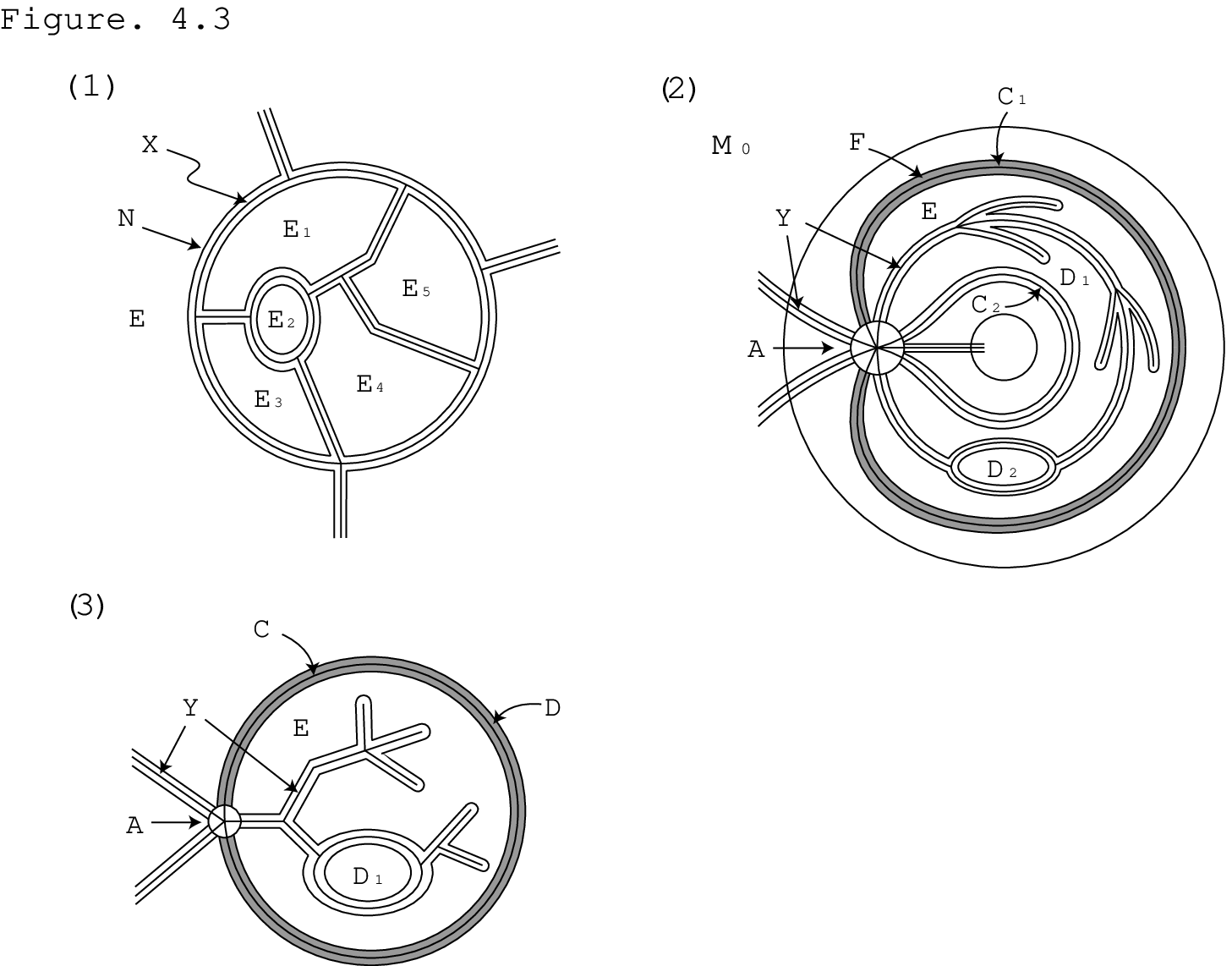} 
\end{center}
\vskip 2mm

\begin{proof}[Proof of Theorem 1.1]
By Lemma 2.1 we may assume that $\partial M = \emptyset$. 
We show that the restriction map $p : {\cal H}(M)_0 \to {\cal E}(X, M)_0$, $p(h) = h|_X$, is a homotopy equivalence.

Let $N$ be a regular neighborhood of $X$ and let $E_i$ ($i = 1, \cdots, m$) and $F_j$ ($j = 1, \cdots, n$) denote the disk or M\"obius band components and the other components of $cl(M \setminus N)$ respectively, 
and let $\displaystyle N_1 = N \cup \left( \cup_{i=1}^{{}\ \, m}\, E_i \right)$. 
By Lemma 4.2 (1) and Proposition 3.1 the following restriction maps are homotopy equivalences:
\[ {\cal E}(N_1, M)_0 \to {\cal E}(N, M)_0 \to {\cal E}(X, M)_0. \] 

Consider the principal bundle 
\[ {\cal H}(M)_0 \cap {\cal H}_{N_1}(M) \hookrightarrow {\cal H}(M)_0 \stackrel{q}{\longrightarrow} {\cal E}(N_1, M)_0, \ \ \ q(h) = h|_X. \]
Since ${i_{N_1}}_\ast\,\pi_1(N_1)$ contains the noncyclic subgroup ${i_X}_\ast\,\pi_1(X)$, 
the submanifold $N_1$ is not a disk, an annulus and a M\"obius band. 
Since $cl(M \setminus N_1) = \cup_{j=1}^{{}\ \ n}\,F_j$ and each $F_j$ is not a disk nor a M\"obius band, 
from Proposition 2.3 it follows that 
${\cal H}(M)_0 \cap {\cal H}_{N_1}(M) = {\cal H}_{N_1}(M)_0 \simeq \ast$ and the map $q$ is a homotopy equivalence. 
Therefore the map $p$ is also a homotopy equivalence. 

Since $\pi_1(M)$ contains the noncyclic subgroup ${i_X}_\ast\,\pi_1(X)$,  
it follows 
that $M \not\cong {\Bbb S}^2$, ${\Bbb P}^2$, ${\Bbb R}^2$, ${\Bbb P}^2 \setminus 1\,pt$, ${\Bbb S}^1 \times (0,1)$.
Therefore, by Proposition 2.2 ${\cal H}(M)_0$ is homotopy equivalent to (a) ${\Bbb T}^2$ if $M \cong {\Bbb T}^2$, (b) ${\Bbb S}^1$ if $M \cong {\Bbb K}^2$, 
and (c) $\ast$ if $M \not\cong {\Bbb T}^2, {\Bbb K}^2$  
\end{proof}

%%%%%%%%%%%%%%%%%%%%%%%%%%%%%%%%%%%%  Section 5  %%%%%%%%%%%%%%%%%%%%%%%%%%%%%%%%%%%%

\section{Embedding spaces of a circle}

The following is the main result of this section, which implies Theorem 1.2.
Suppose $M$ is a connected 2-manifold. 

\renewcommand{\labelenumi}{(\arabic{enumi})}
\begin{theorem} Suppose $C$ is an essential circle in $M$. 
\begin{itemize}
\item[\it (1)] ${\cal E}(C, M)_0 \simeq {\Bbb S}^1$ if $M \not\cong {\Bbb P}^2, {\Bbb T}^2, {\Bbb K}^2$. 
\item[\it (2)] ${\cal E}(C, M)_0 \simeq {\Bbb T}^2$ if $M \cong {\Bbb T}^2$. 
\item[\it (3)] Suppose $M \cong {\Bbb K}^2$. 
\begin{itemize}
\item[\it (i)] ${\cal E}(C, M)_0 \simeq {\Bbb T}^2$ if $C$ is an o.p. nonseparating circle (a meridian).
\item[\it (ii)] ${\cal E}(C, M)_0 \simeq {\Bbb S}^1$ if $C$ is an o.p. separating circle. (an o.p. longitude $=$ a common boundary of two M\"obius bands). 
\item[\it (iii)] ${\cal E}(C, M)_0 \simeq {\Bbb S}^1$ if $C$ is an o.r. circle (an o.r. longitude). 
\end{itemize}
\item[\it (4)] ${\cal E}(C, M)_0 \simeq SO(3)/{\Bbb Z}_2$ if $M \cong {\Bbb P}^2$ 
\end{itemize}
\end{theorem}

Note that $\pi_k(SO(3)) = \pi_k(SO(3)/{\Bbb Z}_2) = \pi_k(S^2) \ (k \geq 3)$, $\pi_2(SO(3)) = \pi_2(SO(3)/{\Bbb Z}_2) = 0$ and 
$\pi_1(SO(3)) \cong {\Bbb Z}_2$, $\pi_1(SO(3)/{\Bbb Z}_2) \cong {\Bbb Z}_4$. 

When a point of $C$ is fixed, we have the following version:

\begin{proposition} 
Suppose $C$ is an essential circle in $M$ and $x \in C$. 
Then ${\cal E}_x(C, M)_0 \simeq \ast$ if $M \not\cong {\Bbb P}^2$ and ${\cal E}_x(C, M)_0 \simeq {\Bbb S}^1$ if $M \cong {\Bbb P}^2$.
\end{proposition} 

\begin{proof} 
By Lemma 2.1 we may assume that $\partial M = \emptyset$. Consider the principal bundle
\[ {\cal G} \equiv {\cal H}_x(M)_0 \cap {\cal H}_C(M) \subset {\cal H}_x(M)_0  \stackrel{p_1}{\to} {\cal E}_x(C, M)_0, \  p_1(h) = h|_C. \]

Below we show that ${\cal G} = {\cal H}_C(M)_0$. Then ${\cal G} \simeq \ast$ by Proposition 2.2, and the map $p_1$ is a homotopy equivalence. 
Since $C$ is essential in $M$, we have $M \not\cong {\Bbb R}^2$, ${\Bbb S}^2$.
Therefore, by Proposition 2.2 ${\cal E}_x(C, M)_0 \simeq {\cal H}_x(M)_0 \simeq \ast$ if $M \not\cong {\Bbb P}^2$ and \  
$\simeq {\Bbb S}^1$ if $M \cong {\Bbb P}^2$.

Let $N$ be a regular neighborhood of $C$ in $M$. 
If $f \in {\cal G}$, then 
(a) $f$ preserves the local orientation at $x$ since $f \simeq_x id_M$, 
(b) $f$ does not interchange the sides of $C$ at $x$ (and every point of $C$) by (a) and $f|_C = id_C$. 
Hence $f$ is isotopic rel $C$ to a $g \in {\cal H}_N(M)$. 
Since $g \simeq_C f \simeq_x id_M$, we have $g \in {\cal H}_x(M)_0 \cap {\cal H}_N(M)$.
If we show that 
\renewcommand{\labelenumi}{$(\ast)$}
\begin{enumerate}
\item \hspace{10mm} ${\cal H}_x(M)_0 \cap {\cal H}_N(M) = {\cal H}_N(M)_0$,
\end{enumerate}
then we have $id_M \simeq_N g \simeq_C f$ and so ${\cal G} = {\cal H}_C(M)_0$. 

The claim ($\ast$) follows from Proposition 2.3$'$ as follows:

(I) If $N$ is an annulus, then we have (a) $x \in N$, (b) $cl(M \setminus N)$ has no disk component since $C$ is essential and 
(c) $cl(M \setminus N)$ has either (i) no M\"obius band component, (ii) exactly one M\"obius band component $L$, and 
(iii) two M\"obius band components $L_1$ and $L_2$. 
In the case (i), ($\ast$) follows from Proposition 2.3$'$.  
In the case (ii), (a) any $h \in {\cal H}_x(M)_0 \cap {\cal H}_N(M)$ is isotopic rel $N$ to a $k \in {\cal H}_{N_1}(M)$ ($N_1 = N \cup L$) 
since ${\cal H}_\partial(L) \simeq \ast$, and 
(b) ${\cal H}_x(M)_0 \cap {\cal H}_{N_1}(M) = {\cal H}_{N_1}(M)$ by Proposition 2.3$'$, so 
(c) $h \in {\cal H}_N(M)_0$ since $k \simeq_N h \simeq_x id_M$ and $k \in {\cal H}_{N_1}(M)_0$ ($k \simeq_{N_1} id_M$).
In the case (iii),  we have ${\cal H}_N(M) = {\cal H}_N(M)_0$ (any $h \in {\cal H}_N(M)$ is isotopic rel $N$ to $id_M$.) since $M = N \cup L_1 \cup L_2$. 
This implies ($\ast$).

(II) If $N$ is a M\"obius band, then $x \in N$ and $L = cl(M \setminus N)$ is connected. 
If $L$ is not a disk nor a M\"obius band, then ($\ast$) follows from Proposition 2.3$'$.
If $L$ is a disk or a M\"obius band, then ${\cal H}_N(M) = {\cal H}_N(M)_0$ and this implies ($\ast$). 
This completes the proof.
\end{proof}

In the proof of Theorem 5.1 we are concerned with the following fiber bundles : 
Suppose $C$ is an essential circle in $M$ and $x \in C$. 
\[ \begin{array}{cccccccccl}
(1) & {\cal F} & \equiv & {\cal E}(C, M)_0 \cap {\cal E}_{x}(C, M) & \subset & {\cal E}(C, M)_0 & \stackrel{p}{\to} & M & : & p(f) = f(x), \\
(2) & {\cal G} & \equiv & {\cal H}_x(M)_0 \cap {\cal H}_C(M) & \subset & {\cal H}_x(M)_0 & \stackrel{p_1}{\to} & {\cal E}_x(C, M)_0 & : & p_1(h) = h|_C \\
(3) & {\cal H} & \equiv & {\cal H}(M)_0 \cap {\cal H}_C(M) & \subset & {\cal H}(M)_0 & \stackrel{p_2}{\to} & {\cal E}(C, M)_0 & : & p_2(h) = h|_C, \\
(4) & {\cal K} & \equiv & {\cal H}(M)_0 \cap {\cal H}_{x}(M) & \subset & {\cal H}(M)_0 & \stackrel{p_3}{\to} & M & : & p_3(h) = h(x). 
\end{array} \]

Suppose $\partial M = \emptyset$ and  
let $\alpha \in \pi_1(M, x)$ be the element represented by $C$ with an orientation.

\begin{lemma} 
(1) If $M \not\cong {\Bbb P}^2$, then 
$\pi_k({\cal E}(C, M)_0) = 0$ $(k \geq 2)$ and $p_\ast : \pi_1({\cal E}(C, M)_0, i_C) \hookrightarrow \pi_1(M, x)$ is a monomorphism. \\
(2) $p_1$ is a homotopy equivalence. \\
(3) $p_{2 \ast} : \pi_k{\cal H}(M)_0 \cong \pi_k{\cal E}(C, M)_0$ $(k \geq 2)$ and 
$p_{2 \ast} : \pi_1{\cal H}(M)_0 \hookrightarrow \pi_1{\cal E}(C, M)_0$ is a monomorphism. \\
(4) (a) $\alpha \in {\rm Im} \, p_{\ast} \subset \pi_1(M, x)$,
(b) $\alpha \beta = \beta \alpha$ for any $\beta \in {\rm Im} \, p_{\ast}$, and 
(c) If $M \not\cong {\Bbb P}^2$, then $\langle \alpha \rangle \cong {\Bbb Z}$. 
\end{lemma}

\begin{proof} 
(1), (3) Since ${\cal F}_0 = {\cal E}_x(C, M)_0 \simeq \ast$ for $M \not\cong {\Bbb P}^2$ (Proposition 5.1) and 
${\cal H}_0 = {\cal H}_C(M)_0 \simeq \ast$, 
the assertions follow from the exact sequences of the fibrations $p$ and $p_2$. 

(2) The assertion has been verified in the proof of Proposition 5.1.

(4)(b) Every map $\phi : ({\Bbb S}^1, \ast) \to ({\cal E}(C, M)_0, i_C)$ induces a map $\Phi : ({\Bbb S} \times C, (\ast, x)) \to (M, x)$. 
Since $p_{\ast}[\phi], \alpha  \in {\rm Im} \, \Phi_{\ast} \subset \pi_1(M, x)$, we have the conclusion. 
\end{proof} 

\begin{proof}[Proof of Theorem 5.1]
By Lemma 2.1 we may assume that $\partial M = \emptyset$. 

(1) Since $M \not\cong {\Bbb P}^2$, by Lemma 5.1 (1), (4) 
it suffices to show that ${\rm Im} \, p_{\ast} \subset \langle \alpha \rangle$ so that ${\rm Im} \, p_{\ast} \cong {\Bbb Z}$. 
Let $\beta \in {\rm Im} \, p_{\ast}$. 
Since $M \not\cong {\Bbb T}^2$, ${\Bbb K}$ and $\alpha \beta = \beta \alpha$ (Lemma 5.1 (4)(b)), 
it follows that $\alpha, \beta \in \langle \delta \rangle$ for some $\delta \in \pi_1(M, x)$ and 
that $\alpha = \delta^k$ and $\beta = \delta^\ell$ for some $k, \ell \in {\Bbb Z}$, $k \neq 0$, so $\alpha^\ell = \beta^k$. 

(i) Suppose $C$ does not bound a M\"obius band. 
Since $M \not\cong {\Bbb K}$ and $C$ is essential, we have $\beta \in \langle \alpha \rangle$ (Lemma 2.2 (v)). 

(ii) Suppose $C$ bounds a M\"obius band $E$. 
Then $C$ is an o.p. circle, and in the proof of Lemma 3.1 (2)(a) Case (iii), 
we have already shown that the loop $pf = f_t(x)$ is o.p. for any class $[f] \in \pi_1({\cal E}(C, M)_0, i_C)$ 
(or for any isotopy $f : [0, 1] \to {\cal E}(C, M)_0$ such that $f_0 = f_1 = i_C$). 
(Note that, for any regular neighborhood $N$ of $C$, using the bundle ${\cal E}(N, M)_0 \to {\cal E}(C, M)_0$, 
we can always find an isotopy $f_t' \in {\cal E}(N, M)_0$ such that $f_0' = i_N$ and $f_t'|_C = f_t$.) 
This observation means that $\beta$ is o.p.  

Let $\gamma \in \pi_1(M, x)$ denote the element 
which is represented by the center circle $A$ of $E$ and satisfies $\gamma^2 = \alpha$. 
Since $M \not\cong {\Bbb K}$ and $A$ does not bound any disk nor M\"obius band, and since $\beta^k = \gamma^{2\ell}$, 
it follows that $\beta \in \langle \gamma \rangle$ (Lemma 2.2 (v)).
Since $\beta$ is o.p. and $\gamma$ is o.r., we have $\beta \in \langle \alpha \rangle$ as required.

(2) If $G$ is a path-connected topological group, then for any point $a \in G$ 
the map $q : ({\cal H}(G)_0, id_G) \to (G, a)$, $q(h) = h(a)$, admits a section $s : (G, a) \to ({\cal H}(G)_0, id_G)$, $s(x)(y) = xa^{-1}y$. 
Hence if $M \cong {\Bbb T}^2$, then $p_{3 \ast} : \pi_1({\cal H}(M)_0, id_M) \to \pi_1(M, x)$ is surjective and 
so is $p_\ast : \pi_1({\cal E}(C, M)_0, i_C) \to \pi_1(M, x)$. 
Hence by Lemma 5.1 (1) $p : {\cal E}(C, M)_0 \to M$ is a homotopy equivalence. 

(3)  
Let $a, b \in \pi_1({\Bbb K})$ denote the classes represented by the meridian $m$ and the o.r. longitude $\ell$ of ${\Bbb K}$ respectively.
Then $\pi_1({\Bbb K}) = \langle a, b \, : \, bab^{-1} = a^{-1} \rangle$ and 
the center of $\pi_1({\Bbb K}) = \langle b^2 \rangle \cong {\Bbb Z}$ \cite[\S3(f)]{Sc}. 

(i) Since $(M, X) \cong ({\Bbb K}, m)$, we may assume that $(M, X) = ({\Bbb K}, m)$ and $\alpha = a$. 
By Lemma 5.1 (1) it suffices to show that ${\rm Im} \, p_{\ast} \cong {\Bbb Z} \oplus {\Bbb Z}$. 
Note that $a \in {\rm Im} \, p_{\ast}$, $b \not\in {\rm Im} \, p_{\ast}$ since $ab \neq ba$ (Lemma 5.1 (4)) and 
$b^2 = p_\ast[f] \in {\rm Im} \, p_{\ast}$, where the loop $f_t \in {\cal E}(m,{\Bbb K})_0$ isotopes $m$ twice along $\ell$.
Therefore ${\rm Im} \, p_{\ast} = \langle a, b^2 \rangle$ (the subgroup of $\pi_1({\Bbb K})$ generated by $a$ and $b^2$).
Since the natural double cover ${\Bbb T}^2 \to {\Bbb K}^2$ corresponds to $\langle a, b^2 \rangle \subset \pi_1({\Bbb K})$, 
we have $\langle a, b^2 \rangle \cong\pi_1({\Bbb T}^2) \cong {\Bbb Z} \oplus {\Bbb Z}$.

(ii) Since $M$ is a union of two M\"obius bands with the common boudary circle $C$, from Lemma 4.2 (1) it follows that 
${\cal E}(C, M)_0 \cong {\cal H}(M)_0 \cong {\Bbb S}^1$.

(iii) Since $(M, X) \cong ({\Bbb K}, \ell)$, we may assume that $(M, X) = ({\Bbb K}, \ell)$ and $\alpha = b$. 
By Lemma 5.1 (1), (4) it suffices to show that ${\rm Im} \, p_{\ast} \subset \langle \alpha \rangle$. 
Given any $\beta = a^r b^s \in {\rm Im} \, p_{\ast}$. 
Since $b \beta = \beta b$, it follows that $a^{2r} = 1$. Since $\pi_1({\Bbb K})$ has no torsion, we have $r = 0$ and $\beta = b^s \in \langle \alpha \rangle$. 

(4) We use the following notations: 
We regard as ${\Bbb R}^3 = {\Bbb C} \times {\Bbb R}$. 
${\Bbb S}^2$ is the unit sphere of ${\Bbb R}^3$. 
 $C_0 = \{ (z, x) \in {\Bbb S}^2 \mid x = 0 \}$ and $N_0 = \{ (z, x) \in {\Bbb S}^2 \mid 0 \leq x \leq 1/2 \}$.
$\pi : {\Bbb S}^2 \to {\Bbb P}^2$ denotes the natural double covering, which identifies antipodal points $(z, x)$ and $(-z, -x)$. 
Since $(M, C) \cong ({\Bbb P}^2, \pi(C_0))$, we may assume that $(M, C) = ({\Bbb P}^2, \pi(C_0))$. 
Then $N = \pi(N_0)$ is a M\"obius band with the center circle $C$. 

We will construct the following diagram: 
\[ \begin{CD}
SO(3) @>{\lambda}>{\simeq}> {\cal H}({\Bbb P}^2)_0 @>{p_4}>{\simeq}> {\cal E}(N,{\Bbb P}^2)_0 @. \\
@V{q_1}VV    @V{q_2}VV   @V{q_3}VV  \\
SO(3)\big/\langle h \rangle @>>{\overline{\lambda}}> {\cal H}({\Bbb P}^2)_0/\langle \overline{h} \rangle @>>{\overline{p_4}}> {\cal E}(N,{\Bbb P}^2)_0/\langle \overline{h}|_N \rangle 
@>{\simeq}>{\overline{p}}> {\cal E}(C, {\Bbb P}^2)_0 \\
\end{CD} \]

Each $f \in SO(3)$ induces a unique $\overline{f} \in {\cal H}({\Bbb P}^2)_0$ with $\pi f = \overline{f} \pi$.
The map $\lambda : SO(3) \to {\cal H}({\Bbb P}^2)_0$, $f \mapsto \overline{f}$, is a homotopy equivalence (Proposition 2.2 (1)(i)).
The restriction map $p_4$ is a homotopy equivalence by Lemma 4.2 (1).

Consider the involution $h \in SO(3)$, $h(z, x) = (-z, x)$. 
By right composition, the group $\langle h \rangle = \{ id_{{\Bbb S}^2}, h \}$ acts on $SO(3)$ 
and $\langle \overline{h} \rangle = \{ id_{{\Bbb P}^2}, \overline{h} \}$ acts on ${\cal H}({\Bbb P}^2)_0$. 
The vertical maps $q_1$ and $q_2$ are the associated quotient maps, which are double coverings.

Since $h(N_0) = N_0$, it follows that $\overline{h}(N) = N$, $\overline{h}|_N \in {\cal H}_C(N)$, $(\overline{h}|_N)^2 = id_N$ and 
$\langle \overline{h}|_N \rangle = \{ id_N, \overline{h}|_N \}$ acts on ${\cal E}(N,{\Bbb P}^2)_0$ by right composition.  
Since $\overline{h}|_N$ interchanges the two local sides of $C$ in $N$, 
by Proposition 3.1 (2) the restriction map $p : {\cal E}(N,{\Bbb P}^2)_0 \to {\cal E}(C, {\Bbb P}^2)_0$ factors 
as the composition of the quotient double covering $q_3$ and the homotopy equivalence $\overline{p}$. 

Since the maps $\lambda$ and $p_4$ are equivariant with respect to these ${\Bbb Z}_2$-actions, 
they induce the associated maps $\overline{\lambda}$ and $\overline{p_4}$. 
Since $\lambda$, $p_4$ are homotopy equivalences and $q_i$'s are covering, 
the maps $\overline{\lambda}$ and $\overline{p_4}$ induce isomorphisms on the $k$-th homotopy groups for $k \geq 2$. 

We show that these maps also induce isomorphisms on $\pi_1$ and so they are homotopy equivalences. 
Since $\pi_1(SO(3)) \cong {\Bbb Z}_2$ and $q_1$ is a double covering, the order $\# \pi_1(SO(3))/\langle h \rangle) = 4$. 
Consider the loop $f_t \in SO(3)$, $f_t(z, x) = (e^{2\pi i t} z, x)$ ($0 \leq t \leq 1$). 
The class $\alpha = [f_t]$ generates $\pi_1(SO(3))$. 
Since $f_{1/2} = h$, the loop $q_1f_t \in SO(3)/\langle h \rangle$ ($0 \leq t \leq 1/2$) induces a class $\beta \in \pi_1(SO(3)/\langle h \rangle)$.
Since $f_{t+1/2} = f_t h$ ($0 \leq t \leq 1/2$), we have $\beta^2 = {q_1}_{\ast} \alpha \neq 1$. 
Therefore ${\text order}\,\beta = 4$ and $\pi_1(SO(3)/\langle h \rangle) = \langle \beta \rangle \cong {\Bbb Z}_4$.

Same argument applies to show that 
$\pi_1({\cal H}({\Bbb P}^2)_0/\langle \overline{h} \rangle) = \langle \overline{\beta} \rangle \cong {\Bbb Z}_4$ and 
$\pi_1({\cal E}(N,{\Bbb P}^2)_0/\langle \overline{h}|_N \rangle) = \langle \overline{\beta}|_N \rangle \cong {\Bbb Z}_4$,
where $\overline{\beta} = [q_2 \overline{f}_t \, (0 \leq t \leq 1/2)]$ and $\overline{\beta}|_N = [q_3 \overline{f}_t|_N \, (0 \leq t \leq 1/2)]$.

Since $\overline{\lambda}_\ast(\beta) = \overline{\beta}$ and $\overline{p_4}_\ast(\overline{\beta}) = \overline{\beta}|_N$, 
it follows that $\overline{\lambda}$ and $\overline{p_4}$ induce isomorphisms on $\pi_1$. 
This completes the proof.
\end{proof} 
  
\begin{proof}[Proof of Theorem 1.2]
When $X$ is a circle, the assertions follow from Theorem 5.1 directly. 
Below we assume that $X$ is not a circle. 
By Lemma 2.1 we may assume that $\partial M = \emptyset$. 

Let $N$ be a regular neighborhood of $X$. 
Since $X$ is not a closed 2-manifold, $N$ is a compact connected 2-manifold with boundary and admits 
a subpolyhedron $Y$ such that $N$ is a regular neighborhood of $Y$ in $M$ and 
$\displaystyle Y = D \cup \left( \cup_{i=1}^m \, C_i \right) \cup \left( \cup_{i=j}^n \, C_j' \right)$
 is a one point union of a disk $D$, essential circles $C_i$ ($i = 1, \cdots m$) ($m \geq 1$) 
and inessential circles $C_j'$ ($j = 1, \cdots n$) ($n \geq 0$). 
Let $Y_1 = D \cup \left( \cup_{i=1}^m \, C_i \right)$.
By Proposition 3.1 and Lemma 4.2 (3) the restriction maps 
\[ {\cal E}(X, M)_0 \longleftarrow {\cal E}(N, M)_0 \longrightarrow {\cal E}(Y, M)_0 \longrightarrow {\cal E}(Y_1, M)_0 \]
are homotopy equivalences.

Note that ${i_X}_\ast \pi_1(X) = {i_N}_\ast \pi_1(N) = {i_Y}_\ast \pi_1(Y) = {i_{Y_1}}_\ast \pi_1(Y_1)$ is a cyclic subgroup of $\pi_1(M)$.
Hence by Lemma 4.1 each pair $(C_k, C_\ell)$ (or $(C_\ell, C_k)$) ($1 \leq k, \ell \leq m$, $k \neq \ell$) satisfies 
one of the conditions of Lemma 4.1 (1), (2) and (3).
By Lemma 4.2 (2) there exists a $k$ ($1 \leq k \leq m$) such that the restriction map 
\[ {\cal E}(Y_1, M)_0 \longrightarrow {\cal E}(D \cup C_k, M)_0 \]
is a homotopy equivalence. 

Let $N_1$ be a regular neighborhood of $D \cup C_k$. 
Then $N_1$ is an annulus or a M\"obius band, which is a regular neighborhood of $C_k$.
We set $A = C_k$ when $N_1$ is an annulus and $A = \partial N$ when $N_1$ is a M\"obius band. 
By Proposition 3.1 and Lemma 4.2 (1) the restriction maps
\[  {\cal E}(D \cup C_k, M)_0 \longleftarrow {\cal E}(N_1, M)_0 \longrightarrow {\cal E}(A, M)_0 \]
are homotopy equivalences.

We apply Theorem 5.1 to the circle $A$.
The statements (1), (2) follow from Theorem 5.1 (1), (2) directly. 

(3) Suppose $M \cong {\Bbb K}$. 

(i) Suppose $X$ is contained in an annulus $N_0$ which does not separate $M$. 
We may assume that $X \subset Int\,N_0$ and $N \subset N_0$.  
Since $C_k$ is essential and $C_k \subset N_0$ it follows that $C_k$ is an o.p. nonseparating circle. 
Hence $A = C_k$ and ${\cal E}(A, M)_0 \cong {\Bbb T}^2$. 

(ii) Suppose $C_k$ is an o.p. nonseparating circle of $M$. 
Each inessential circle $C_j'$ bounds a disk $E_j$. 
Since each $C_i$ is essential, every disk $E_j$ does not intersect $\cup_i\, C_i$ except the wedge point of $Y$.
Since $C_k$ is o.p., the choice of $C_k$ means that 
each $C_i$ is also o.p. and each pair $(C_i, C_k)$ satisfies the condition of Lemma 4.1 (1). 
Hence we can find an annulus neighborhood $N_0$ of $C_k$ with $Y \subset Int\,N_0$.
Since $C_k$ is nonseparating, so is $N_0$. 
Since $N$ is a regular neighborhood of $Y$, we can isotope $N_0$ so that $N \subset N_0$.   

This observation means that if $X$ does not satisfy the condition (3)(i), then $C_k$ is either (a) o.p. and separating or (b) o.r.
In the case (a), $A = C_k$, and in the case (b), $N_1$ is a M\"obius band and $A = \partial N_1$.
In each case, $A$ is an o.p. separating circle and ${\cal E}(A, M)_0 \simeq {\Bbb S}^1$. 

(4) Suppose $M \cong {\Bbb P}^2$. 
Since $A$ is an o.p. circle, $M$ is the union of a disk and a M\"obius band with a common boundary $A$. 
By Lemma 4.2 (1) and Proposition 2.2 (1)(i) ${\cal E}(A, M)_0 \simeq {\cal H}(M)_0 \simeq SO(3)$.
This completes the proof of Theorem 1.2. 
\end{proof}

%%%%%%%%%%%%%%%%%%%%%%%%%%%%%%%%%%%%%%%%%% Subsection 6 %%%%%%%%%%%%%%%%%%%%%%%%%%%%%%%%%%%%%%%%%%

\section{Embedding spaces of an arc and a disk} 

%%%%%%%%%%%%%%%%%%%%%%%%%%%%%%%%%%%%%%%%%% Subsection 6.1 %%%%%%%%%%%%%%%%%%%%%%%%%%%%%%%%%%%%%%%%%%

\subsection{Main statements} \mbox{}
\vskip 1mm
Suppose $M$ is a $2$-manifold and $X$ is a compact connected polyhedron ($\neq$ 1pt) in $M$ with a distinguished point $x \in X$. 
In this section we identify the fiber homotopy (f.h.) type of the projection $p : {\cal E}(X, M)_0 \to M$, $p(f) = f(x)$ in the case where $X \simeq \ast$ in $M$.

We choose a smooth structure and a Riemannian metric of $M$ and consider the unit circle bundle $q : S(TM) \to M$ of the tangent bundle $q : TM \to M$. 
Let $\pi : \tilde{M} \to M$ denote the orientation double cover of $M$, which has a natural orientation and the Riemannian metric induced from $M$. 
Let $\tilde{q} : S(T\tilde{M}) \to \tilde{M}$ denote the associated unit circle bundle of $\tilde{M}$. 
The terminology ``fiber homotopy equivalence (or equivalent)'' is abbreviated as f.h.e.

\begin{theorem} Suppose $X \simeq \ast$ in $M$ and $\partial M = \emptyset$. 
Then $p : {\cal E}(X, M)_0 \to M$ is f.h.e over $M$ to \\
(i) $q : S(TM) \to M$ if $X$ is an arc or $M$ is orientable, \\
(ii) $\pi\tilde{q} : S(T\tilde{M}) \to M$ if $X$ is not an arc and $M$ is nonorientable. 
\end{theorem} 

Theorem 1.3 follows from Theorem 6.1, Lemma 2.1 and the fact that $S(TM) \simeq S(TInt\,M)$.
Since $X$ has a disk neighborhood in $M$, Theorem 6.1 is reduced to the following more technical propositions: 
Suppose $D$ is an oriented disk and $X$ is a compact connected polyhedron ($\neq$ 1pt) in $Int\,D$ 
with a distinguished point $x \in X$. 
Consider the subspace
\[ {\cal E}^\ast(X, M) = \{ f \in {\cal E}(X, M) \mid \ \mbox{$f$ admits an extension $\overline{f} \in {\cal E}(D, M)$}\}, \]
and the projection $p : {\cal E}^\ast(X, M) \to M$, $p(f) = f(x)$. 
When $M$ is oriented, consider the subspaces 
\[ {\cal E}^{\pm}(X, M) = \{ f \in {\cal E}(X, M) \mid \ \mbox{$f$ admits an o.p./o.r. extension $\overline{f} \in {\cal E}(D, M)$}\}. \]
Since $\tilde{M}$ has a natural orientation, this definition applies to spaces of embeddings into $\tilde{M}$.
Let $\tilde{p} : {\cal E}(X, \tilde{M}) \to \tilde{M}$, $\tilde{p}(f) = f(x)$ denote the projection.

\begin{proposition} 
(1) Suppose $X$ is not an arc. \\
(i) When $M$ is oriented, the projection $p : {\cal E}^{\pm}(X, M) \to M$ is f.h.e. to $q : S(TM) \to M$ over $M$. \\
(ii) When $M$ is nonorientable, $p : {\cal E}^\ast(X, M) \to M$ is f.h.e. to $\pi \tilde{q} : S(T\tilde{M}) \to M$ over $M$. \\
(2) When $X$ is an arc, the projection $p : {\cal E}(X, M) \to M$ is f.h.e. to $q : S(TM) \to M$ over $M$. 
\end{proposition}

When $X$ is an arc and $x$ is an interior point of $X$, we can introduce a ${\Bbb Z}_2$-action: 
Let $I = [-1, 1]$ and let $p : {\cal E}(I, M) \to M$, $p(f) = f(0)$ denote the projection.
The group ${\Bbb Z}_2 = \{ \pm 1 \}$ admits 
a f.p. action on ${\cal E}(I, M)$ by $(\varepsilon \cdot f)(t) = f(\varepsilon t)$ ($\varepsilon \in {\Bbb Z}_2$, $t \in I$), and 
a f.p. action on $S(TM)$ by $\varepsilon \cdot v = \varepsilon v$ ($\varepsilon \in {\Bbb Z}_2$, $v \in S(T_xM)$, $x \in M$).

\begin{proposition}
The projection $p : {\cal E}(I, M) \to M$ is ${\Bbb Z}_2$-equivariant f.h.e. to $q : S(TM) \to M$ over $M$. 
\end{proposition}

Propositions 6.1 and 6.2 will be verified in the subsequent subsections. 

%%%%%%%%%%%%%%%%%%%%%%%%%%%%%%%%%%%%%%%%%% Subsection 6.2 %%%%%%%%%%%%%%%%%%%%%%%%%%%%%%%%%%%%%%%%%%

\subsection{Extension Lemma} \mbox{} 
\vskip 1mm 
Let $D(1)$ ($O(1)$) denote the closed (open) unit disk in ${\Bbb R}^2$ and suppose $X$ is a compact connected polyhedron in $O(1)$. 
In this subsection we apply the conformal mapping theorem in the complex function theory 
so as to construct a canonical extension map $\Phi : {\cal E}^\ast(X, O(1)) \to {\cal H}(D(1))$ and show the naturality and symmetry properties of $\Phi$.
The case where $X$ is a tree has been treated in \cite{Ya2}. 

%%%%%%%%%%%%%%%%%%%%%%%%%%%%%%%%%%%%%%%%%% Subsection 6.2.1 %%%%%%%%%%%%%%%%%%%%%%%%%%%%%%%%%%%%%%%%%%

\subsubsection{Canonical parametrizations} \mbox{} 
\vskip 1mm 
Suppose $Y$ is a compact 1-dim polyhedron.
Let $R(Y)$ denote the set of points of $Y$ which have a neighborhood homeomorphic to ${\Bbb R}$, and set $V(Y) = Y \setminus R(Y)$. 
Each point of $V(Y)$ is called a vertex of $Y$ and the closure of each component of $R(Y)$ in $Y$ is called an edge of $Y$. 
Therefore, an edge $e$ is an arc or a simple closed curve: 
in the former case the end points of $e$ are vertices and in the latter case $e$ contains at most one vertex. 
By $E(Y)$ we denote the set of edges of $Y$. 
Note that $V(Y)$ and $E(Y)$ are topological invariants of $Y$. 
%For any edge $e \in E(Y)$, the open subset $\stackrel{\circ}{e} = e \setminus V(Y)$ is called the open edge of $e$.  
An oriented edge $e$ of $Y$ means an edge of $Y$ with a distinguished orientation.
By $e^{-1}$ we denote the same edge with the opposite orientation. 

Suppose $X$ is a compact connected polyhedron ($\neq$ 1 pt) topologically embedded in $O(1)$.   
It follows that $X$ is a subpolyhedron with respect to some triangulation of $O(1)$ and that 
$O(1) \setminus X$ is a disjoint union of an open annulus $U^X$ and a finite number of open disks. 
Let $\Lambda(X)$ and $\Lambda_0(X)$ denote the set of all components and the subset of open disk components of $O(1) \setminus X$. 

For each $U \in \Lambda(X)$, ${\rm Fr}\,U \ (= {\rm Fr}_{O(1)}U)$ is a compact connected 1-dim polyhedron. 
For the annulus (respectively each disk) component $U \in \Lambda(X)$ let ${\cal E}(U)$ denote the set of oriented edges $e$ of ${\rm Fr}\,U$ such that 
the right (respectively left) hand side of $e$ lies in $U$. 
Then ${\cal E}(U)$ admits a unique cyclic ordering ${\cal E}(U) = \{ e_U(1), \cdots, e_U(n_U) \}$ 
($e_U(n_U+1) = e_U(1)$) such that 
\begin{itemize}
\item[$(\ast)$] $e_U(j)$ and $e_U(j+1)$ are adjacent when they are seen from $U$, and have compatible orientations for $j = 1, \cdots, n_U$. 
\end{itemize}

If we move on these edges in this order, we obtain a loop $\ell_U$ which moves on ${\rm Fr}\,U$ in the ``counterclockwise'' orientation. 
As a normalization data, for each $U \in \Lambda_0(X)$ 
we choose an ordered set $a_U = (x_U, y_U, z_U)$ of three distinct points lying on the loop $\ell_U$ in the positive order, 
while on $C(1)$ we take the ordered set $a_0 = (-i, 1, i)$ (invariant under $\eta$).

The conformal mapping theorem yields a canonical parametrization of each $U \in \Lambda(X)$. 
Based on the boundary behaviours of these conformal mappings, we obtain the next lemmas:

\begin{lemma} 
(1) For the annulus component $U = U^X$, there exists a unique $r = r_X \in (0, 1)$ and 
a unique o.p. map $g = g_X : A(r, 1) \to cl\,U \subset D(1)$ 
such that $g$ maps ${\rm Int}\,A(r, 1)$ conformally onto $U$ and $g(1) = 1$. 
Furthermore, $g$ satisfies the following conditions: 
(a) $g$ maps $C(1)$ homeomorphically onto $C(1)$,  
(b) $g(C(r)) = {\rm Fr}U$ and $g$ satisfies the condition $(\#)_U$ on $C(r)$.

(2) For each disk component $U \in \Lambda_0(X)$, 
there exists a unique o.p. map $g = g_{(U, a_U)} : D(1) \to cl\,U \subset D(1)$ 
such that $g$ maps $O(1)$ conformally onto $U$ and $g(a_0) = a_U$. 
Furthermore, $g$ satisfies the following conditions: 
$g(C(1)) = {\rm Fr}U$ and $g$ satisfies the condition $(\#)_U$ on $C(1)$.
\end{lemma} 
\noindent Here, the condition $(\#)_U$ on $C(r)$ is stated as follows: 
\begin{itemize}
\item[$(\#)_U$] There exists a unique collection of points $\{ u_U(1), \cdots, u_U(n_U) \}$  
lying on $C(r)$ in counterclockwise order such that 
$g$ maps each positively oriented circular arc $\overline{u_U(j)u_U(j+1)}$ onto the oriented edge $e_U(j)$ in o.p. way and 
maps $\displaystyle Int\,\left[ \overline{u_U(j)u_U(j+1)} \right]$ homeomorphically onto $e_U(j) \setminus V({\rm Fr}U)$.
(Here $u_U(n_U +1) = u_U(1)$, and when $n_U = 1$, we mean that $\overline{u_U(1)u_U(1)} = C(r)$.)
\end{itemize}
\vskip 2mm

For $0 < r < 1$ we define a radial map $\lambda_r : A(1/2, 1) \to A(r, 1)$ by $\lambda_r(x) = (2(1 - r)(|x| - 1) + 1)x/|x|$. 
We set $h_X = g_X \lambda_{r_X} \in {\cal C}(A(1/2, 1), D(1))$.  

%%%%%%%%%%%%%%%%%%%%%%%%%%%%%%%%%%%%%%%%%% Subsection 6.2.2 %%%%%%%%%%%%%%%%%%%%%%%%%%%%%%%%%%%%%%%%%%

\subsubsection{Canonical extensions} \mbox{} 
\vskip 1mm 
Suppose $(X, a)$ and $(Y, b)$ are two compact connected polyhedra ($\neq$ 1 pt) in $O(1)$ 
with normalization data $a = \{ a_U \}_{U \in \Lambda_0(X)}$ and $b = \{ b_V \}_{V \in \Lambda_0(Y)}$.
In the case where $X$ is not an arc, if $f : X \to Y$ is any homeomorphism which admits an extension $\overline{f} \in {\cal H}(D(1))$, 
then the sign $\delta(f) = \pm$ is defined by $\overline{f} \in {\cal H}^{\delta(f)}(D(1))$, which depends only on $f$.  
In the case where $X$ and $Y$ are arcs, we consider any pair $(f, \delta)$ of a homeomorphism $f : X \to Y$ and $\delta = \pm$. 
By abuse of notation, $(f, \delta)$ is simply denoted by $f$, and $\delta$ by $\delta(f)$. 
In this setting we will construct a canonical extension $\Phi_{a,b}(f) \in {\cal H}^{\delta(f)}(D(1))$ of $f$.

By the choice of $\delta$, we can find an extension $\overline{f} \in {\cal H}^\delta(D(1))$ of $f$. 
The subsequent arguments do not depend on the choice of such an extension $\overline{f}$.  
The statement $A$/$B$ mean that $A$ holds for $\delta = +$ and $B$ holds for $\delta = -$.

(1) For each disk component $U \in \Lambda_0(X)$, consider the corresponding disk component $U_f = \overline{f}(U) \in \Lambda_0(Y)$ 
($U_f$ is independent of the choice of $\overline{f}$). 
Lemma 6.1 (2) provides with two maps $g_{(U, a_U)}$ and $g_{(U_f, b_{U_f})}$. 

(2) For the annulus component $U_X$, consider the corresponding annulus component $U_Y$. 
Lemma 6.1 (1) provides with two data $(r_X, g_X, h_X)$ and $(r_Y, g_Y, h_Y)$. 
For the notational compatibility, for $U = U_X$, let $U_f = U_Y$, $g_{U} = g_X$ and $g_{U_f} = g_Y$.

For any $U \in \Lambda(X)$, it follows that $\overline{f} : (U, {\rm Fr}\,U) \cong (U_f, {\rm Fr}\,U_f)$ is an o.p./o.r. homeomorphism, and that
if $\{ e(1), \cdots, e(n) \}$ is the cyclic ordering of ${\cal E}(U)$, 
then $\{ f(e(1)), \cdots, f(e(n)) \}$ represents the positive/negative cyclic ordering of ${\cal E}(U_f)$. 
In particular, it also follows that $f(\ell_U) = (\ell_{U_f})^\delta$. 
Thus, by reversing the orientation and order for $\delta = -$, the condition ``$(\#)_{U_f}$ on $C(r_f)$'' ($r_f = 1$ or $r_Y$) can be restated as follows:
\begin{itemize}
\item[$(\#\#)_{U_f}$] There exists a unique collection of points $\{ v(1), \cdots, v(n) \}$  
lying on $C(r_f)$ in counterclockwise/clockwise order such that 
$g_f$ maps each oriented circular arc $\overline{v(j)v(j+1)}$ onto the oriented edge $f(e(j))$ in o.p. way and 
maps $Int\,\left[\overline{v(j)v(j+1)}\right]$ homeomorphically onto $f(e(j) \setminus V({\rm Fr}U))$.
(As before, $v(n +1) = v(1)$, and when $n = 1$, we mean that $\overline{v(1)v(1)} = C(r)^{\delta}$.)
\end{itemize}
\vskip 2mm 
(1) For each $U \in \Lambda_0(X)$, compare two maps $fg_U$, $g_{U_f} : C(1) \to {\rm Fr}U_f$. 
By the conditions $(\#)_U$ and $(\#\#)_{U_f}$, we obtain a unique map $\theta_U(f) \in {\mathcal H}^\delta(C(1))$ such that $g_{U_f} \theta_U(f) = f g_U$. 
Extend $\theta_U(f)$ conically to 
\[ \Theta_U(f) \in {\mathcal H}^\delta(D(1)): \ \ \Theta_U(f)(sz) = s \theta_U(f)(z) \ \ (z \in C(1), 0 \leq s \leq 1). \]
There exists a unique homeomorphism $\phi_U(f) : cl(U) \cong cl(U_f)$ which satisfies $g_{U_f} \Theta_U(f) = \phi_U(f) g_U$.  
Then $\phi_U(f)$ is an extension of $f : {\rm Fr}U \cong {\rm Fr}U_f$. 

\noindent (2) Compare two maps $fh_X$, $h_Y : C(1/2) \to {\rm Fr}U_Y$. 
By the conditions $(\#)_{U_X}$ and $(\#\#)_{U_Y}$, we obtain a unique map $\theta_X(f) \in {\mathcal H}^\delta(C(1/2))$ such that $h_Y\theta_X(f) = fh_X$. 
This definition can be also applied when $X$ is an arc. 
Extend $\theta_X(f)$ radially to 
\[ \Theta_X(f) \in {\mathcal H}^\delta(A(1/2, 1)): \ \ \Theta_X(f)(sz) = s\theta_X(f)(z/2) \ \ (z \in C(1), 1/2 \leq s \leq 1). \] 
There exists a unique homeomorphism $\phi_X(f) : cl(U_X) \cong cl(U_Y)$ which satisfies $h_Y\Theta_X(f) = \phi_X(f) h_X$.  
Then $\phi_X(f)$ is an extension of $f : {\rm Fr}U_X \cong {\rm Fr}U_Y$. 

Finally we define $\Phi(f) \in {\cal H}(D(1))$ by $\Phi(f) = f$ on $X$, $\Phi(f) = \phi_U(f)$ on $cl(U)$ ($U \in \Lambda_0(X)$) 
and $\Phi(f) = \phi_X(f)$ on $cl(U_X)$.  

%%%%%%%%%%%%%%%%%%%%%%%%%%%%%%%%%%%%%%%%%% Subsection 6.2.3 %%%%%%%%%%%%%%%%%%%%%%%%%%%%%%%%%%%%%%%%%%

\subsubsection{Extension map} \mbox{}
\vskip 1mm
Suppose $(X, a)$ is a compact connected polyhedron ($\neq$ 1 pt) in $O(1)$ with a normalization data $a = \{ a_U \}_{U \in \Lambda_0(X)}$.
First suppose $X$ is not an arc. We keep the notations given in \S6.1, so that 
\renewcommand{\labelenumi}{(\roman{enumi})}
\begin{enumerate}
\item ${\cal E}^{\pm}(X, O(1)) = \{ f \in {\cal E}(X, O(1)) \mid \mbox{ $f$ admits an extension $\overline{f} \in {\cal H}^{\pm}(D(1))$}\}$,
\item ${\cal E}^\ast(X, O(1)) = {\cal E}^+(X, O(1)) \cup {\cal E}^-(X, O(1))$ and ${\cal E}^+(X, O(1)) \cap {\cal E}^-(X, O(1)) = \emptyset$,
\item ${\cal E}^+(X, O(1)) = {\cal E}(X, O(1))_0$ and ${\cal E}^{\pm}(X, O(1)) = \{ \eta f \mid f \in {\cal E}^{\mp}(X, O(1))\}$,
\end{enumerate}
where $\eta : {\Bbb R}^2 \cong {\Bbb R}^2$ denotes the reflection $\eta(x, y) = (x, -y)$. 
In the statement (iii), the former follows from the bundle $p : {\cal H}(D(1))_0 \to {\cal E}(X, D(1))_0$, $p(h) = h|_X$. 

For any $f \in {\cal E}^\delta(X, O(1))$, 
the image $f(X)$ is a compact connected polyhedron in $O(1)$, 
to which we can assign a normalization data $a_f = \{ (a_f)_{U_f} \}_{U \in \Lambda_0(X)}$ defined by the condition:   
\[ (a_f)_{U_f} = f(a_U)^\delta \mbox{($\equiv (f(x), f(y), f(z))$ for $\delta = +$ and $\equiv (f(z), f(y), f(x))$ for $\delta = -$).} \]
Since $f(\ell_U) = (\ell_{U_f})^\delta$, $(a_f)_{U_f}$ lies on $\ell_{U_f}$ in the positive order. 

By \S6.2.2 we obtain (i) the parametrization: 
$(r_{f(X)}, g_{f(X)}, h_{f(X)})$ for the annulus component $U_{f(X)}$ and $g_{U_f}$ for each disk component $U_f \in \Lambda_0(f(X))$, and 
(ii) the canonical extension $\Phi_{(X,a)}(f) = \Phi_{a, a(f)}(f) \in {\cal H}^\delta(D(1))$ of $f$. 

When $X$ is an arc, we set ${\cal E}^\ast(X, O(1)) = {\cal E}(X, O(1)) \times \{ \pm \}$ 
(note that ${\cal E}(X, O(1)) = {\cal E}^{\pm}(X, O(1))$ in the usual sense).
We identify ${\cal E}(X, O(1))$ with ${\cal E}^+(X, O(1))$. 

By \S6.2.2, for any $(f, \delta) \in {\cal E}^\ast(X, O(1))$, we obtain 
(i) the parametrization: $(r_{f(X)}, g_{f(X)}, h_{f(X)})$ for the annulus component $U_{f(X)}$, and 
(ii) the canonical extension $\Phi_X(f, \delta) \in {\cal H}^\delta(D(1))$ of $f$. 
In the subsequent statements the notation $(X, a)$ simply means $X$ when $X$ is an arc. 
The next assertions follow from the same argument as in \cite[Lemma 2.3, Proposition 2.1]{Ya2}. 

\begin{lemma}
The next correspondence is continuous: 
\[ {\cal E}^\ast(X, O(1)) \ni f \mapsto (r_{f(X)}, h_{f(X)}, \{ g_{U_f} \}_{U \in \Lambda_0(X)}) 
\in (0, 1) \times {\cal C}(A(1/2, 1), D(1)) \times \prod_{U \in \Lambda_0(X)} {\cal C}(D(1), D(1)). \] 
\end{lemma}

\begin{proposition}
The correspondence $\Phi_{(X,a)} : {\cal E}^\ast(X, O(1)) \to {\cal H}(D(1))$ is continuous. 
\end{proposition} 

The extension map $\Phi$ has the following naturality and symmetry properties: 
For notations: As usual we identify ${\Bbb C}$ with ${\Bbb R}^2$. 
For each $z \in C(1)$, let $\theta_z : {\Bbb C} \cong {\Bbb C}$ denote the rotation $\theta_z(w) = z \cdot w$. 
Let $\eta : {\Bbb R}^2 \cong {\Bbb R}^2$ denote the reflection $\eta(x, y) = (x, -y)$. 
Let $\eta_+ = id$ and $\eta_- = \eta$. 
The restriction of $\gamma \in O_2$ to $E = C(r), D(r), O(r)$ is denoted by the same symbol $\gamma$. 

The group $O_2$ acts on ${\cal H}(E)$ by left and right composition and on ${\cal E}(X, E)$ by left composition.    
These conventions yield the natural embeddings $\theta : C(1) \cong SO_2 \subset {\cal H}(E)_0$ and $O_2 \subset {\cal H}(E)$.  
We regard as $SO_2 \subset {\cal E}(X, O(1))_0$ by $SO_2 \ni \theta \mapsto \theta|_X \in {\cal E}(X, O(1))_0$. 

If $X = [-1/2, 1/2] \subset D(1)$, then $\theta_{-1}$ preserves $X$ and  
the group ${\Bbb Z}_2 = \{ \pm 1 \}$ admits right actions on ${\cal E}^\ast(X, E)$ and ${\cal H}(E)$ 
by $(f, \delta) \cdot \varepsilon = (f \theta_{\varepsilon}|_X, \delta)$ and $h \cdot \varepsilon = h \theta_{\varepsilon}$. 
These actions restrict to the $SO_2$-actions and ${\Bbb Z}_2$-actions on ${\cal H}^+(E)$ and ${\cal E}^+(X, E)$.

\begin{proposition} Suppose $X$ is not an arc. \\
(1) $\Phi_{(X,a)}(gf) = \Phi_{(f(X), a(f))}(g) \Phi_{(X,a)}(f)$ for any $f \in {\cal E}^\ast(X, O(1))$ and $g \in {\cal E}^\ast(f(X), O(1))$. \\
(2) $\Phi_{(X,a)}(\gamma|_X) = \gamma$ for any $\gamma \in O_2$. \\
(3) $\Phi_{(X,a)} : {\cal E}^\ast(X, O(1)) \to {\cal H}(D(1))$ is left $O_2$-equivariant. 
\end{proposition}

\begin{proposition} Suppose $X$ is an arc. \\
(1) $\Phi_X(gf, \varepsilon \delta) = \Phi_{f(X)}(g, \varepsilon) \Phi_X(f, \delta)$ 
for any $(f, \delta) \in {\cal E}^\ast(X, O(1))$ and $(g, \varepsilon) \in {\cal E}^\ast(f(X), O(1))$. \\
(2) $\Phi_X(\gamma|_X, \delta(\gamma)) = \gamma$ for any $\gamma \in O_2$. \\
(3) (i)(a) $\Phi_X : {\cal E}^\ast(X, O(1)) \to {\cal H}(D(1))$ is left $O_2$-equivariant. \\
\hspace*{9.5mm} (b) $\Phi_{\eta(X)}(\eta f \eta |_{\eta(X)}, \delta) = \eta \Phi_X(f, \delta) \eta$ for any $(f, \delta) \in {\cal E}^\ast(X, O(1))$. \\
\hspace*{4mm} (ii) If $X = [-1/2, 1/2]$, then \\
\hspace*{9.5mm} (a) $\Phi_X : {\cal E}^\ast(X, O(1)) \to {\cal H}(D(1))$ is right ${\Bbb Z}_2$-equivariant, and (b) $\Phi_X(i_X, -) = \eta$. 
\end{proposition}

\begin{proof}[Proof of Propositions 6.4 and 6.5]
(1) Since $\delta(gf) = \delta(g) \delta(f)$, $U_{gf} = (U_f)_g$ ($U \in \Lambda(X)$) and $a_{gf} = (a_f)_g$, it follows that 
$\theta_U(gf) = \theta_{U_f}(g) \theta_U(f)$, $\Theta_U(gf) = \Theta_{U_f}(g) \Theta_U(f)$ and $\phi_U(gf) = \phi_{U_f}(g) \phi_U(f)$.  
This implies that $\Phi_{(X,a)}(gf) = \Phi_{(f(X), a(f))}(g) \Phi_{(X,a)}(f)$ 
(or $\Phi_X(gf, \varepsilon \delta) = \Phi_{f(X)}(g, \varepsilon) \Phi_X(f, \delta)$).

(2) It suffices to show that 
(a) $\Phi_{(X,a)}(\theta_z|_X) = \theta_z$, (b) $\Phi_{(X,a)}(\eta|_X) = \eta$ 
(and (a) $\Phi_X(\theta_z|_X, +) = \theta_z$, (b) $\Phi_X(\eta|_X, -) = \eta$ when $X$ is an arc).

(a) For each $U \in \Lambda_0(X)$, the uniqueness part of Lemma 6.1 (1) implies that $h_{U_{\theta_z}} = \theta_z h_U$. 
This means that $\Theta_U(\theta_z|_X) = id$ and $\phi_U(\theta_z|_X) = \theta_z$.

For $U = U_X$, let $w \in C(1)$ be the unique point such that $\theta_z g_X \theta_z^{-1}(w) = 1$. 
The uniqueness part of Lemma 6.1 (2) implies that $g_{\theta_z(X)} = \theta_z g_X \theta_{z^{-1}w}$ and $h_{\theta_z(X)} = \theta_z h_X \theta_{z^{-1}w}$.  
This means that $\Theta_X(\theta_z|_X) = \theta_{w^{-1}z}$ and $\phi_X(\theta_z|_X) = \theta_z$. 
Thus we have $\Phi(\theta_z|_X) = \theta_z$.

(b) For each $U \in \Lambda(X)$, the uniqueness part of Lemma 6.1 implies that $\eta h_X = h_{\eta(X)} \eta$. 
This means that $\theta_U(\eta|_X) = \eta$, $\Theta_U(\eta|_X) = \eta$ and $\phi_U(\eta|_X) = \eta$. 
Thus we have $\Phi(\eta|_X) = \eta$.

(3) The assertions follow form (1) and (2). 
Note that $(f \theta_{-1}|_X, \delta) = (f, \delta) (\theta_{-1}|_X, +)$ and $(\eta f \eta, \delta) = (\eta, -)(f, \delta)(\eta, -)$ when $X = [-1/2, 1/2]$. 
\end{proof}

%%%%%%%%%%%%%%%%%%%%%%%%%%%%%%%%%%%%%%%%%% Subsection 6.3 %%%%%%%%%%%%%%%%%%%%%%%%%%%%%%%%%%%%%%%%%%

\subsection{Deformation Lemma} \mbox{} 
\vskip 1mm 
Suppose $X$ is a compact connected polyhedron ($\neq$ 1 pt) in $O(1)$ with a normalization data. 

%%%%%%%%%%%%%%%%%%%%%%%%%%%%%%%%%%%%%%%%%% Subsection 6.3.1 %%%%%%%%%%%%%%%%%%%%%%%%%%%%%%%%%%%%%%%%%%

\subsubsection{Deformation of ${\cal E}^+(X, O(1))$ onto a circle} \mbox{} 
\vskip 1mm 
In this subsection we use the extension maps $\Phi_X$  
to construct an $O_2$-equivariant strong deformation retraction ($O_2$-s.d.r.) of the embedding space ${\cal E}^\ast(X, O(1))$ onto $O_2$.
The space $O_2$ is embedded into ${\cal H}(E)$ ($E = C(1)$, $D(1)$) by the restriction $\gamma \mapsto \gamma|_E$ and 
into ${\cal E}^\ast(X, O(1))$ by $\gamma \mapsto \gamma|_X$ (and $\gamma \mapsto (\gamma|_X, \delta(\gamma))$ when $X$ is an arc).

We need some auxiliary homotopies $G_t$, $A_t$ and $H_t$ ($0 \leq t \leq 1$):  \\
(i) ${\cal H}(C(1))$ has a natural s.d.r. $G_t$ onto $O_2$ defined by 
\[ G_t : {\cal H}(C(1)) \to {\cal H}(C(1)), \ \  
\ G_t(g)(e^{i \theta}) = g(1) \,\exp\,[i\delta(g)((1-t) \tau(\theta) + t \theta)], \]
where 
$\delta(g) = \pm$ is defined according to $g \in {\cal H}^{\pm}(C(1))$ and 
$\tau : [0, 2\pi] \to [0, 2\pi]$ is a unique map such that $\tau(0) = 0$ and $g(e^{i \theta}) = g(1) \exp\,[i\delta(g)\tau(\theta)]$. \\
(ii) The cone extension map $c : {\cal H}(C(1)) \to {\cal H}(D(1))$, $c(g)(sx) = sg(x)$ \ ($x \in C(1)$, $0 \leq s \leq 1$), 
is a section of the restriction map $p : {\cal H}(D(1)) \to {\cal H}(C(1))$, $p(h) = h|_{C(1)}$. \\
(iii) The Alexander trick yields a s.d.r. $A_t$ of ${\cal H}(D(1))$ onto ${\rm Im}\,c$ defined by 
\[ A_t : {\cal H}(D(1)) \to {\cal H}(D(1)) : \ id \simeq cp, \ \ 
A_t(h)(x) = \begin{cases} 
|x|h(x/|x|) & (s \leq |x| \leq 1, x \neq 0) \\
sh(x/s) & (0 < |x| \leq s) \\
0 & (x = 0, t = 1)
\end{cases} (\mbox{where } s = 1 - t). 
 \]
(iv) Combining these homotopies we obtain a s.d.r. $H_t$ of ${\cal H}(D(1))$ onto $O_2$: 
\[ H_t : {\cal H}(D(1)) \to {\cal H}(D(1)): \ \ H_t = 
\begin{cases} 
A_{2t} & (0 \leq t \leq 1/2) \\
cG_{2t-1}p & (1/2 \leq t \leq 1) 
\end{cases}. 
\] 

The maps $c$, $p$ and $A_t$ are easily seen to be left and right $O_2$-equivariant. 

\begin{lemma} 
$G_t$ and $H_t$ are left $O_2$ and right $\eta$-equivariant. 
\end{lemma}

\begin{proof}
$G_t$ : (i) Let $\gamma = \theta_z$ ($z = e^{i \lambda}$). 
Then 
\[ (\gamma g)(e^{i\theta}) = (\gamma g)(1) e^{i \delta(g) \tau(\theta)} \ \mbox{ and } \ 
G_t(\gamma g)(e^{i \theta}) = e^{i\, \lambda} g(1) \exp[i \delta(g) ( (1 - t)\tau(\theta) + t \theta )] = \gamma G_t(g)(e^{i \theta}). \]

\noindent (ii) Since $\eta(e^{i \theta}) = e^{-i \theta} = e^{i(2\pi - \theta)}$ and $\delta(\eta g) = - \delta(g)$, it follows that 
\begin{eqnarray*}
(a) \hspace{15mm} (\eta g)(e^{i \theta}) &=& (\eta g)(1) \exp[i \delta(\eta g) \tau(\theta)], \\
G_t(\eta g)(e^{i \theta}) &=& (\eta g)(1) \exp[i \delta(\eta g)((1-t)\tau(\theta) + t \theta)] = (\eta G_t(g))(e^{i \theta}), \\
(b) \hspace{15mm} (g\eta)(e^{i \theta}) &=& g(1) e^{i\delta(g)\tau(2\pi - \theta)} = (g \eta)(1) e^{i\delta(g \eta)(2\pi - \tau(2\pi - \theta))}, \\ 
G_t(g \eta)(e^{i \theta}) &=& (g \eta)(1) \exp[i \delta(g \eta)( (1-t)(2\pi - \tau(2\pi - \theta)) + t \theta)] \\
&=& g(1) \exp[i \delta(g)( (1-t)\tau(2\pi - \theta) + t(2\pi - \theta))] \\
&=& G_t(g)(e^{i (2\pi - \theta)}) = G_t(g)(\eta(e^{i \theta})). 
\end{eqnarray*}
\vskip -10mm
\end{proof}

The required $O_2$-s.d.r. $F_t$ of ${\cal E}^\ast(X, O(1))$ onto $O_2$ is defined by 
\begin{itemize}
\item[(1)] $X$ is not an arc : \ \ \ $\displaystyle F_t : {\cal E}^\ast(X, O(1)) \to {\cal E}^\ast(X, O(1)), \ F_t(f) = H_t(\Phi_X(f))|_X$ 
\item[(2)] $X$ is an arc : $\displaystyle F_t : {\cal E}^\ast(X, O(1)) \to {\cal E}^\ast(X, O(1)), \ F_t(f, \delta) = (H_t(\Phi_X(f, \delta))|_X, \delta)$ 
\end{itemize}

The $O_2$-equivariance follows from Propositions 6.4, 6.5 and Lemma 6.3. 
When $0 \in X$, we can consider the subspace ${\cal E}^{\pm}(X,0; O(1),0) = \{ f \in {\cal E}^{\pm}(X; O(1)) \mid f(0) = 0 \}$. 

\begin{lemma}
If $0 \in X$, then $F_t({\cal E}^{\pm}(X,0; O(1),0)) \subset {\cal E}^{\pm}(X,0; O(1),0)$. 
\end{lemma}

The circle $C(1)$ is embedded into ${\cal E}^+(X, O(1))$ by $z \mapsto \theta_z|_X$, 
which corresponds to the embedding $SO_2 \subset {\cal E}^+(X, O(1))$.
The $O_2$-s.d.r. $F_t$ restricts to the $SO_2$-s.d.r. of ${\cal E}^+(X, O(1))$ onto $C(1)$. 
In the case where $X$ is an arc, The $SO_2$ action on ${\cal E}(X, O(1))$ ($= {\cal E}^+(X, O(1))$) extends to the $O_2$-action by left composition
($\gamma f = \gamma \circ f$ or $\gamma (f, +) = (\gamma f, +)$ even if $\delta(\gamma) = -$). 
(This action should be distinguished from the $O_2$-action on ${\cal E}^\ast(X, O(1))$.)

\begin{lemma}
When $X$ is an arc $[a, b]$ ($-1 < a < b < 1$), the s.d.r. $F_t : {\cal E}(X, O(1)) \to {\cal E}(X, O(1))$ is left $O_2$-equivariant.
\end{lemma}

\begin{proof}
Since $F_t$ is $SO_2$-equivariant, it suffices to show that $F_t(\eta f) = \eta F_t(f)$ ($f \in {\cal E}(X, O(1))$).
Since $\Phi_X(\eta f, +) = \Phi_X(\eta f \eta|_X, +) = \eta \Phi_X(f, +)\eta$, from Lemma 6.3 it follows that 
\vskip 1mm
$\displaystyle F_t(\eta f) = H_t(\Phi_X(\eta f, +))|_X = H_t(\eta \Phi_X(f, +)\eta)|_X = \eta H_t(\Phi_X(f, +)) \eta|_X = \eta F_t(f)$. 
\end{proof}

When $X$ is the arc $I = [-1/2, 1/2]$, we can modify the construction of $F_t$ in order that $F_t$ is right ${\Bbb Z}_2$-equivariant: 
Let $J = \{ \pm 1 \} \subset C(1)$. \\
(i) A s.d.r. $\nu_t$ of ${\cal E}(J, C(1))$ onto $C(1)$ is defined by 
\[ \nu_t : {\cal E}(J, C(1)) \to {\cal E}(J, C(1)): \ \nu_t(\alpha)(\pm 1) = \alpha(\pm 1) \, \exp[\mp it(\pi/2 - \theta(\alpha))], \] 
where $\theta(\alpha) \in (0, 2\pi)$ is defined by $\displaystyle \alpha(-1) = \alpha(1) e^{2i\theta(\alpha)}$. \\
(ii) The cone extension map $c' : {\cal E}(J, C(1)) \to {\cal E}([-1, 1], D(1))$ is defined by 
$c'(\alpha) (\pm s) = s \alpha(\pm 1)$ \ ($s \in [0,1]$). \\
(iii) Let $q : {\cal H}(D(1)) \to {\cal E}(J, C(1))$ denote the restriction map. 

The modified s.d.r. $F_t$ of ${\cal E}(I, O(1))$ onto $C(1)$ is defined by 
\[ F_t : {\cal E}(I, O(1)) \to {\cal E}(I, O(1)), \ 
F_t(f) = \begin{cases}
\left. A_{2t}(\Phi_I(f, +)) \right|_I & (0 \leq t \leq 1/2) \\[1mm]
\left.c'(\nu_t q(\Phi_I(f, +))) \right|_I & (1/2 \leq t \leq 1) 
\end{cases} \]

\begin{lemma}
(i) $F_t$ is left $O_2$, right ${\Bbb Z}_2$-equivariant. \\
(ii) $F_t({\cal E}(I,0; O(1),0)) \subset {\cal E}(I,0; O(1),0)$.
\end{lemma}

\begin{proof}
(i) The map $A_t$ is left and right $O_2$-equivariant, and the maps $c'$, $q$, $\nu_t$ are easily seen to be left $O_2$, right ${\Bbb Z}_2$-equivariant. 
Since $\Phi_I(\ , +)$ is left $SO_2$, right ${\Bbb Z}_2$-equivariant (Proposition 6.5 (3)), so is $F_t$. 
The left $\eta$-equivariance of $F_t$ follows from Proposition 6.5 (3)(i)(b) and the direct observation:
\[ F_t(\eta f) = \begin{cases}
A_{2t}(\eta \Phi_I(f, +) \eta)|_I = \eta A_{2t}(\Phi_I(f, +)) \eta |_I = \eta A_{2t}(\Phi_I(f, +))|_I \\[1mm]
c'(\nu_t q(\eta \Phi_I(f, +) \eta))|_I = c'(\nu_t \eta q(\Phi_I(f, +)))|_I = \eta c'(\nu_t q(\Phi_I(f, +))) |_I
\end{cases} \]
\end{proof}

%%%%%%%%%%%%%%%%%%%%%%%%%%%%%%%%%%%%%%%%%% Subsection 6.3.2 %%%%%%%%%%%%%%%%%%%%%%%%%%%%%%%%%%%%%%%%%%

\subsubsection{Oriented plane case} \mbox{} 
\vskip 1mm
Suppose $V$ is an oriented 2-dim vector space with an inner product. 
For $\varepsilon \in (0, \infty]$, 
let $D_V(\varepsilon) = \{ v \in V \mid |v| \leq \varepsilon \}$, $O_V(\varepsilon) = \{ v \in V \mid |v| < \varepsilon \}$ and 
$C_V(\varepsilon) = \{ v \in V \mid |v| = \varepsilon \}$.
Consider the subspace 
\[ {\cal E}^\pm(X, O_V(\varepsilon)) 
= \{ f \in {\cal E}(X, O_V(\varepsilon)) \mid \mbox{$f$ extends to an o.p./o.r. homeomorphism $\overline{f} : D(1) \cong D_V(\varepsilon)$} \}. \]
The circle $C_V(1)$ is embedded into ${\cal E}^+(X, O_V(\varepsilon))$ by $v \mapsto \varepsilon\alpha_v|_X$, 
where $\alpha_v : {\Bbb R}^2 \cong V$ is the unique o.p. linear isometry with $\alpha_v(1) = v$. 
The group $SO(V)$ acts on ${\cal E}^+(X, O_V(\varepsilon))$ by left composition ($\gamma f = \gamma \circ f$).
The $SO_2$-s.d.r. $F_t$ of ${\cal E}^+(X, O(1))$ onto $C(1)$ induces a $SO(V)$-s.d.r. $\phi_t^V$ of ${\cal E}^+(X, O_V(1))$ onto $C_V(1)$:
\[ \phi_t^V : {\cal E}^+(X, O_V(1)) \to {\cal E}^+(X, O_V(1)), \ \phi_t^V(f) = \alpha F_t(\alpha^{-1}f), \]
where $\alpha : {\Bbb R}^2 \cong V$ is any o.p. linear isometry. 
Since $F_t$ is $SO_2$-equivariant, $\phi_t^V$ does not depend on the choice of $\alpha$. 
When $0 \in X$, we have $\phi_t^V({\cal E}^+(X,0; O_V(1),0)) \subset {\cal E}^+(X,0; O_V(1),0)$ (Lemma 6.4). 

The group $SO(V)$ is canonically isomorphic to $SO_2$ under 
the isomorphism $\chi_V : SO_2 \cong SO(V)$, $\chi_V(\gamma) = \alpha \gamma \alpha^{-1}$.
Here $\alpha : {\Bbb R}^2 \cong V$ is any o.p. linear isometry, and $\chi_V$ does not depend on the choice of $\alpha$. 
Thus $SO_2$ acts on ${\cal E}^+(X, O_V(\varepsilon))$ canonically and $\phi_t^V$ is a $SO_2$-s.d.r.

When $X = [a, b]$, we can include the unoriented case.
Suppose $V$ is a 2-dim (unoriented) vector space with an inner product. 
The circle $C_V(1)$ is embedded into ${\cal E}(X, O_V(\varepsilon))$ by $v \mapsto \varepsilon\alpha_v|_X$, 
where $\alpha_v : {\Bbb R}^2 \cong V$ is any one of two linear isometries with $\alpha_v(1) = v$. 
The s.d.r. $F_t : {\cal E}(X, O(1)) \to {\cal E}(X, O(1))$ is $O_2$-equivariant (Lemma 6.5) and it induces 
a $O(V)$-s.d.r. $\phi_t^V$ of ${\cal E}(X, O_V(1))$ onto $C_V(1)$:
\[ \phi_t^V : {\cal E}(X, O_V(1)) \to {\cal E}(X, O_V(1)), \ \phi_t^V(f) = \alpha F_t(\alpha^{-1}f), \]
where $\alpha : {\Bbb R}^2 \cong V$ is any linear isometry and $\phi_t^V$ does not depend on the choice of $\alpha$. 
When $0 \in X$, we have $\phi_t^V({\cal E}(X,0; O_V(1),0)) \subset {\cal E}(X,0; O_V(1),0)$. 
In the unoriented case there is no canonical isomorphism $SO_2 \cong SO(V)$. 

When $X$ is the arc $I = [-1/2, 1/2] \subset O(1)$, 
we can use the modified s.d.r. $F_t$ so that $\phi_t^V$ is left $O(V)$ and right ${\Bbb Z}_2$-equivariant.

%%%%%%%%%%%%%%%%%%%%%%%%%%%%%%%%%%%%%%%%%% Subsection 6.3.3 %%%%%%%%%%%%%%%%%%%%%%%%%%%%%%%%%%%%%%%%%%

\subsubsection{Oriented plane bundle case} \mbox{} 
\vskip 1mm
The purpose of this subsection is to construct a fiberwise deformation of 
the fiberwise embedding space ${\cal E}^+_\pi(X, O_E(\varepsilon))$ onto the unit circle bundle $S(E)$ for any oriented plane bundle $p : E \to B$. 

Suppose $\pi : E \to B$ is an oriented 2-dim vector bundle with an inner product. 
Let $V_b = \pi^{-1}(b)$ for $b \in B$. 
For any map $\varepsilon : B \to (0, \infty)$, we can associate the subspaces: 
\vskip 2mm
\[
\begin{array}{l}
O_E(\varepsilon) = \cup_{b \in B}\, O_{V_b}(\varepsilon(b)) \subset E, \ \ \ \ 
{\cal E}_{\pi}^{(+)}(X, O_E(\varepsilon)) = \cup_{b \in B}\,{\cal E}^{(+)}(X, O_{V_b}(\varepsilon(b)) \subset {\cal C}(X, E), \\[2mm]
{\cal E}_{\pi}^{(+)}(X,0; O_E(\varepsilon),0) = \cup_{b \in B}\,{\cal E}^{(+)}(X,0; O_{V_b}(\varepsilon(b)),0) 
\subset {\cal E}_{\pi}^{(+)}(X, O_E(\varepsilon)) \ \text{ when } 0 \in X, 
\end{array}
\]
\vskip 2mm
\noindent and the projection: \ \ 
$\displaystyle p : {\cal E}_{\pi}^{(+)}(X, O_E(\varepsilon)) \to B$, $p(f) = b$ \ \ ($f \in {\cal E}(X, O_{V_b}(\varepsilon(b))$). 

Let $\pi : S(E) \to B$ denote the unit circle bundle of the bundle $\pi : E \to B$.  
The $SO_2$-actions on $O_{V_b}(1)$ and $C_{V_b}(1)$ ($b \in B$) induce 
a f.p. $SO_2$-actions on $O_E(\varepsilon)$, $S(E)$ and ${\cal E}_{\pi}^+(X, O_E(\varepsilon))$. 
The $SO_2$-embeddings $C_{V_b}(1) \to {\cal E}^+(X, O_{V_b}(\varepsilon))$ induce a f.p. $SO_2$-embedding $S(E) \to {\cal E}_{\pi}^+(X, O_E(\varepsilon))$. 
The $SO_2$-s.d.r.'s $\Phi_t^b = \phi_t^{V_b}$ of ${\cal E}^+(X, O_{V_b}(1))$ onto $C_{V_b}(1)$ ($b \in B$) induce  
a f.p. $SO_2$-s.d.r. of ${\cal E}_{\pi}^+(X, O_E(1))$ onto $S(E)$:
\[ \Phi_t : {\cal E}_{\pi}^+(X, O_E(1)) \to {\cal E}_{\pi}^+(X, O_E(1)), \ \ \Phi_t(f) = \Phi_t^b(f) \ (f \in {\cal E}^+(X, O_{V_b}(1))), \]
When $0 \in X$, we have $\Phi_t({\cal E}_{\pi}^+(X,0; O_E(1),0)) \subset {\cal E}_{\pi}^+(X,0; O_E(1),0))$, 
where $0 \subset O_E(1)$ is the image of the zero-section of $E$. 

We define a fiberwise radial homeomorphism $k_\varepsilon : O_E(1) \cong O_E(\varepsilon)$ by $k_\varepsilon(sv) = \varepsilon(b)v$ ($v \in O_{V_b}(1)$). 
Then $k$ is a f.p. $SO_2$-homeomorphism and induces 
a f.p. $SO_2$-homeomorphism $K_\varepsilon : {\cal E}_{\pi}^+(X, O_E(1)) \cong {\cal E}_{\pi}^+(X, O_E(\varepsilon))$, $K_\varepsilon(f) = k_\varepsilon \circ f$.
When $0 \in X$, we have $K({\cal E}_{\pi}^+(X,0; O_E(1),0)) = {\cal E}_{\pi}^+(X,0; O_E(\varepsilon),0)$. 

Finally, the required f.p. $SO_2$-s.d.r. $\Phi_t^\varepsilon$ of ${\cal E}_{\pi}^+(X, O_E(\varepsilon))$ onto $S(E)$ is defined by 
\[ \Phi_t^\varepsilon : {\cal E}_{\pi}^+(X, O_E(\varepsilon)) \to {\cal E}_{\pi}^+(X, O_E(\varepsilon)), 
\ \ \ \Phi_t^\varepsilon = K_\varepsilon \Phi_t K_\varepsilon^{-1}. \]
When $0 \in X$, we have $\Phi_t^\varepsilon({\cal E}_{\pi}^+(X,0; O_E(1),0)) \subset {\cal E}_{\pi}^+(X,0; O_E(\varepsilon),0)$.

When $X = [a, b]$, it follows that ${\cal E}_{\pi}^+(X, O_E(\varepsilon)) = {\cal E}_{\pi}(X, O_E(\varepsilon))$ and 
for any 2-dim vector bundle $\pi : E \to B$ with an inner product, we obtain 
the f.p. s.d.r. $\Phi_t^\varepsilon$ of ${\cal E}_{\pi}(X, O_E(\varepsilon))$ onto $S(E)$. 

When $X$ is the arc $I = [-1/2, 1/2]$, we can obtain a ${\Bbb Z}_2$-equivariant version.  
Suppose $\pi : E \to B$ is a 2-dim vector bundle with an inner product. 
Then ${\cal E}_{\pi}^+(I, O_E(\varepsilon)) = {\cal E}_{\pi}(I, O_E(\varepsilon))$ and 
the right-${\Bbb Z}_2$-s.d.r.'s $\Phi_t^b = \phi_t^{V_b}$ of ${\cal E}(I, O_{V_b}(1))$ onto $C_{V_b}(1)$ ($b \in B$) induces 
a f.p. right-${\Bbb Z}_2$-s.d.r. of ${\cal E}_{\pi}(I, O_E(1))$ onto $S(E)$
\[ \Phi_t : {\cal E}_{\pi}(I, O_E(1)) \to {\cal E}_{\pi}(I, O_E(1)), \ \ \Phi_t(f) = \Phi_t^b(f) \ (f \in {\cal E}(E, O_{V_b}(1))). \]
We have $\Phi_t({\cal E}_{\pi}(I,0; O_E(1),0)) \subset {\cal E}_{\pi}(I,0; O_E(1),0)$. 

The f.p. homeomorphism $K : {\cal E}_{\pi}(I, O_V(1)) \cong {\cal E}_{\pi}(I, O_V(\varepsilon))$ is 
right-${\Bbb Z}_2$-equivariant and $K({\cal E}_{\pi}(I, 0; O_E(1),0)) = {\cal E}_{\pi}(I, 0; O_E(\varepsilon),0)$. 

The required f.p. right ${\Bbb Z}_2$-s.d.r. $\Phi_t^\varepsilon$ of ${\cal E}_{\pi}(I, O_E(\varepsilon))$ onto $S(E)$ is defined by 
\[ \Phi_t^\varepsilon : {\cal E}_{\pi}(I, O_E(\varepsilon)) \to {\cal E}_{\pi}(I, O_E(\varepsilon)), \ \ \Phi_t^\varepsilon = K \Phi_t K^{-1}. \]
We have $\Phi_t^\varepsilon({\cal E}_{\pi}(I,0; O_E(\varepsilon),0)) \subset {\cal E}_{\pi}(I,0; O_E(\varepsilon),0)$.

Replacing $D(1)$ by any oriented disk, we have the following conclusions (the $-$-case is reduced to the $+$-case by reversing the orientation of $D$):

\begin{proposition}
(1) Suppose $D$ is an oriented disk and $X$ is a compact conencted polyhedron ($\neq$ 1pt) in $Int\,D$ with a distinguished point $x_0 \in X$.
Then for any oriented 2-dim vector bundle $\pi : E \to B$ with an inner product, 
there exists a f.p. $SO_2$-s.d.r of ${\cal E}_{\pi}^\pm(X, O_E(\varepsilon))$ (and ${\cal E}_{\pi}^\pm(X,x_0; O_E(\varepsilon),0)$) onto $S(E)$.

(2) Suppose $X$ is an arc and $x_0$ is any point of $X$. 
Then for any 2-dim vector bundle $\pi : E \to B$ with an inner product, 
there exists a f.p. s.d.r of ${\cal E}_{\pi}(X, O_E(\varepsilon))$ (and ${\cal E}_{\pi}(X,x_0; O_E(\varepsilon),0)$) onto $S(E)$.

(3) Let $I = [-1, 1]$ and choose 0 as the base point. Then for any 2-dim vector bundle $\pi : E \to B$ with an inner product, 
there exists a f.p. right ${\Bbb Z}_2$-s.d.r of 
${\cal E}_{\pi}(I, O_E(\varepsilon))$ (and ${\cal E}_{\pi}(I,0; O_E(\varepsilon),0)$) onto $S(E)$.
\end{proposition}

%%%%%%%%%%%%%%%%%%%%%%%%%%%%%%%%%%%%%%%%%% Subsection 6.4 %%%%%%%%%%%%%%%%%%%%%%%%%%%%%%%%%%%%%%%%%%

\subsection{Proof of Theorem 6.1} \mbox{} 

%%%%%%%%%%%%%%%%%%%%%%%%%%%%%%%%%%%%%%%%%% Subsection 6.4.1 %%%%%%%%%%%%%%%%%%%%%%%%%%%%%%%%%%%%%%%%%%
\subsubsection{Spaces of small embeddings} \mbox{} 
\vskip 1mm
Suppose $M$ is a connected $2$-manifold with $\partial M = \emptyset$.
We choose a smooth structure and a Riemannian metric on $M$, which  induces the path-length metric $d$ on $M$. 
The tangent bundle $q : TM \to M$ is a 2-dim vector bundle with an inner product and it is oriented when $M$ is oriented. 
We apply the argument in 6.3.3 to this vector bundle $TM$. 

For $x \in M$ and $r > 0$, let $U_x(r) = \{ y \in M \mid d(x, y) < r \}$ and $O_x(r) = O_{T_x M}(r)$ ($= \{ v \in T_x M \mid |v| < r \}$). 
For any map $\varepsilon : M \to (0, \infty)$, 
let $U_M(\varepsilon) = \cup_{x \in M} \, \{ x \} \times U_x(\varepsilon(x)) \subset M \times M$, 
while $O_{TM}(\varepsilon) = \cup_{x \in M} \, O_x(\varepsilon(x))$ by definition. 
If the map $\varepsilon$ is sufficiently small, then 
at each $x \in M$ the exponential map, ${\rm exp}_x$, maps $O_x(\varepsilon(x))$ diffeomorphically onto $U_x(\varepsilon(x))$ 
(${\rm exp}_x$ is o.p. when $M$ is oriented) \cite[Theorem 1.6]{Cha}. 
Since ${\rm exp}_x$ is smooth in $x \in M$, we obtain a f.p. diffeomorphism over $M$: 
\[ exp : O_{TM}(\varepsilon) \to U_M(\varepsilon), \ \ {\rm exp}(v) = (x, {\rm exp}_x(v)) \ \ (v \in O_x(\varepsilon(x))). \]

In order to connect the space ${\cal E}^\ast(X, M)$ with the fiberwise embedding space ${\cal E}_q^\ast(X, x_0; O_{TM}(\varepsilon), 0)$, 
we introduce spaces of small embeddings. 

Suppose ${\cal U} = \{ U_\lambda \}$ is a cover of $M$ by {\em open disks} and $\delta : M \to (0, \infty)$ is a map with $\delta \leq \varepsilon$. 
For any pointed space $(Y, y_0)$ we set 
\[ \mbox{${\cal E}_{\cal U}(Y, M) = \{ f \in {\cal E}(Y, M) \mid f(Y) \subset U_\lambda$ for some $\lambda\}$, \  

${\cal E}_\delta(Y, M) = \{ f \in {\cal E}(Y, M) \mid f(Y) \subset U_{f(y_0)}(\delta(f(y_0)))\}$.} \] 

Suppose $D$ is an oriented disk and $(X, x_0)$ is a pointed compact connected polyhedron ($\neq$ 1pt) in $Int\,D$ 
($x_0$ is also regarded as a base point of $D$). We consider the following subspace of ${\cal E}^\ast(X, M)$: 
\[ \mbox{${\cal E}^\ast_{\cal U}(X, M) = \{ f \in {\cal E}(X, M) \mid f$ admits an extension $\overline{f} \in {\cal E}_{\cal U}(D, M)\}$.}\]  
When $M$ is oriented, 
we have the subspaces ${\cal E}^{\pm}_{\cal U}(D, M) = {\cal E}_{\cal U}(D, M) \cap {\cal E}^{\pm}(D, M)$ and 
\[ \mbox{${\cal E}^{\pm}_{\cal U}(X, M) = \{ f \in {\cal E}(X, M) \mid f$ admits an extension $\overline{f} \in {\cal E}^{\pm}_{\cal U}(D, M)\}$.}\]
The subspaces ${\cal E}^\ast_\delta(X, M)$ and ${\cal E}^{\pm}_\delta(X, M)$ are defined similarly. 
When $M$ is oriented and $X$ is not an arc, 
we have ${\cal E}^\ast_\delta(X, M) = {\cal E}^+_\delta(X, M) \cup {\cal E}^-_\delta(X, M)$ (a disjoint union). 
Note that if $(M, X) = ({\Bbb S}^2$, a circle) and ${\cal U}$ consists of small open disks, then 
${\cal E}^+_{\cal U}(X, M) \neq {\cal E}^+(X, M) \cap {\cal E}_{\cal U}(X, M)$. 
%
%We note that (i) ${\cal E}^\ast_{\cal U}(X, M) = {\cal E}^\ast(X, M) \cap {\cal E}_{\cal U}(X, M)$ for any ${\cal U}$ iff 
%$M \not\cong {\Bbb S}^2$ or $X$ is contractible or $X$ is a circle and 
%(ii) when $M$ is oriented, ${\cal E}^+_{\cal U}(X, M) = {\cal E}^+(X, M) \cap {\cal E}_{\cal U}(X, M)$ for any ${\cal U}$ iff 
%$M \not\cong {\Bbb S}^2$ or $X$ is contractible. 

The projection $p : {\cal E}(X, M) \to M$, $p(f) = f(x_0)$ restricts to the projections on these subspaces. 

\begin{lemma} 
The f.p. diffeomorphism $exp$ induces a f.p. homeomorphism over $M$ (which preserves the $\pm$-parts when $M$ is oriented): 
\[ Exp : {\cal E}_q^\ast(X, x_0; O_{TM}(\varepsilon), 0) \cong {\cal E}_\varepsilon^\ast(X, M), \ \ \ 
Exp(f) = {\rm exp}_x \circ f \ \ \ (f \in {\cal E}^\ast(X,x_0; O_x(\varepsilon(x)),0))). \]
\end{lemma} 

%Suppose ${\cal U} = \{ U_\lambda \}$ is a cover of $M$ by open disks and $\delta : M \to (0, \infty)$ is a map with $\delta \leq \varepsilon$. 

Let $I$ denote the interval $[-1, 1]$ with the base point $0$. 

\begin{proposition} 
(i)(a) When $M$ is oriented, 
the inclusions ${\cal E}_\delta^\pm(X, M) \subset {\cal E}^\pm(X, M)$ and 
${\cal E}_{\cal U}^\pm(X, M) \subset {\cal E}^\pm(X, M)$ are f.h.e.'s over $M$. \\
(b) When $M$ is nonorientable, the inclusions ${\cal E}_\delta^\ast(X, M) \subset {\cal E}^\ast(X, M)$ and 
${\cal E}_{\cal U}^\ast(X, M) \subset {\cal E}^\ast(X, M)$ are f.h.e.'s over $M$. \\
(ii) The inclusions ${\cal E}_\delta(I, M) \subset {\cal E}(I, M)$ and ${\cal E}_{\cal U}(I, M) \subset {\cal E}(I, M)$ are ${\Bbb Z}_2$-f.h.e.'s over $M$.
\end{proposition} 

First we prove the following assertions. 

\begin{lemma} 
(1) (i) The inclusion ${\cal E}_\delta^\ast(X, M) \subset {\cal E}_\varepsilon^\ast(X, M)$ is a f.h.e. over $M$ 
(which preserves the $\pm$-parts when $M$ is oriented). \\
(ii) The inclusion ${\cal E}_\delta(I, M) \subset {\cal E}_\varepsilon(I, M)$ is a ${\Bbb Z}_2$-f.h.e. over $M$. \\ 
(2) Suppose ${\cal U}_\delta \equiv \{ U_x(\delta(x)) \}_{x \in M}$ refines ${\cal U}$ 
(i.e., each $x \in M$ admits a $\lambda$ with $U_x(\delta(x)) \subset U_\lambda$). \\
(i)(a$_\pm$) When $M$ is oriented, the inclusion ${\cal E}_\delta^\pm(X, M) \subset {\cal E}_{\cal U}^\pm(X, M)$ is a f.h.e. over $M$. \\
(b) When $M$ is nonorientable, the inclusion ${\cal E}_\delta^\ast(X, M) \subset {\cal E}_{\cal U}^\ast(X, M)$ is a f.h.e. over $M$. \\
(ii) The inclusions ${\cal E}_\delta(I, M) \subset {\cal E}_{\cal U}(I, M)$ is a ${\Bbb Z}_2$-f.h.e. over $M$. 
\end{lemma} 

\begin{proof}
(1) Using the radial shrinking of $O_{TM}(\varepsilon)$ onto $O_{TM}(\delta)$, it is seen that the inclusion 
\[ {\cal E}_q^\ast(X, x_0; O_{TM}(\delta), 0) \subset {\cal E}_q^\ast(X, x_0; O_{TM}(\varepsilon), 0) \]
is a (${\Bbb Z}_2$) f.h.e.
Since $Exp : {\cal E}_q^\ast(X, x_0; O_{TM}(\varepsilon), 0) \cong {\cal E}_\varepsilon^\ast(X, M)$ is a f.p. (${\Bbb Z}_2$) homeomorphism and 
$Exp({\cal E}_q^\ast(X, x_0; O_{TM}(\delta), 0)) = {\cal E}_\delta^\ast(X, M)$, we have the conclusion.

(2)(i) 
First we consider the case where $X$ is not an arc. 
We may assume that $D$ is a subdisk of $M$ and that the inclusion $D \subset M$ is o.p. when $M$ is oriented.
We treat the cases (${\rm a}_\pm$) and (b) simultaneously and use the supersprict $\#$ to denote $\pm$ in the case (${\rm a}_\pm$) and $\ast$ in the case (b). 

We show that the restriction map $p : {\cal E}^\#(D, M) \to {\cal E}^\#(X, M)$ has a section $s$ in each case of (${\rm a}_\pm$) and (b). 
If $N_0$ is a regular neighborhood of $X$ in $D$, then $cl(D \setminus N_0)$ is a finite disjoint union of an annulus $A$ and closed disks in $Int\,D$. 
Since $N = N_0 \cup A$ is a regular neighborhood of $X$ in $M$, by Proposition 3.1 (1) and Lemma 4.2 (1)  the restriction maps 
${\cal E}(D, M)_0 \to {\cal E}(N, M)_0 \to {\cal E}(X, M)_0$ are homotopy equivalences. 
Thus the restriction map $p_0 : {\cal E}(D, M)_0 \to {\cal E}(X, M)_0$ is also a homotopy equivalence. 
Since $p_0$ is a locally trivial bundle (Corollary 2.1 (ii)), it has a section $s_0$. 
In the cases (${\rm a}_+$) and (b) it follows that $p = p_0$ (i.e., ${\cal E}(D, M)_0 = {\cal E}^\#(D, M)$, ${\cal E}(X, M)_0 = {\cal E}^\#(X, M)$), 
so we have done. 
The (${\rm a}_-$) case is deduced by taking a $\eta \in {\cal H}^-(D)$ and applying (${\rm a}_+$) to $(D, \eta(X))$.  

Next we show that the restriction map 
$p : {\cal E}^\#_{\cal U}(D, M) \to {\cal E}^\#_{\cal U}(X, M)$ 
has a section $F$ such that $F({\cal E}^\#_{\delta/2}(X, M)) \subset {\cal E}^\#_\delta(D, M)$. 
Since $D \setminus X$ consists of a half-open annulus component $A$ with $\partial A = \partial D$ 
and open disk components $V_i$'s, by shrinking $A$ towards $X$ 
we can find a s.d.r. $\phi_t$ ($0 \leq t \leq 1$) of $D$ onto $\tilde{X} = X \cup (\cup_i V_i)$ such that $\phi_t \in {\cal E}^+(D, D)$ ($0 \leq t < 1$) 
and $\phi_{t_1}(D) \supset \phi_{t_2}(D)$ ($0 \leq t_1 \leq t_2 \leq 1$).

We construct maps $\mu_1 : {\cal E}_{\cal U}^\#(X, M) \to [0, 1)$ and $\mu_2 : {\cal E}_\delta^\#(X, M) \to [0, 1)$ such that 
$s(f)\phi_{\mu_1(f)} \in {\cal E}_{\cal U}^\#(D, M)$ ($f \in {\cal E}_{\cal U}^\#(X, M)$) and 
$s(f)\phi_{\mu_2(f)} \in {\cal E}_\delta^\#(D, M)$ ($f \in {\cal E}_\delta^\#(X, M)$). 
The maps $\mu_1$ is constructed as follows: Given $f \in {\cal E}_{\cal U}^\#(X, M)$, there exists a $\overline{f} \in {\cal E}_{\cal U}^\#(D, M)$
such that $\overline{f}|X = f$. 
By definition $\overline{f}(D) \subset U_\lambda$ for some $\lambda = \lambda_f$. 
Comparing $s(f) \in {\cal E}^\#(D, M)$ and $\overline{f}$, 
we conclude that $s(f)(V_i) = \overline{f}(V_i)$ for each open disk component $V_i$ of $D \setminus X$, so that $s(f)(\tilde{X}) \subset U_{\lambda_f}$.
Thus, if $t_f \in [0, 1)$ is sufficiently close to 1, then $s(f)\phi_{t_f}(D) \subset U_{\lambda_f}$. 
If ${\cal V}_f$ is a sufficiently small neighborhood of $f$ in ${\cal E}_{\cal U}^\#(X, M)$, 
then for each $g \in {\cal V}_f$ we have $s(g)\phi_{t_f}(D) \subset U_{\lambda_f}$.
Choose a locally finite open covering $\{ {\cal W}_f \}$ of ${\cal E}_{\cal U}^\#(X, M)$ 
such that ${\cal W}_f \subset {\cal V}_f$ for each $f$ (${\cal W}_f$ may be empty), and then 
construct a map $\mu_1 : {\cal E}_{\cal U}^\#(X, M) \to [0, 1)$ such that $\mu_1|_{{\cal W}_f} \geq t_f$ for each $f$. 
Then $\mu_1$ satisfies the required condition. 
(For $\mu_2$, replace ${\cal U}$ by $\delta$ and $U_{\lambda_f}$ by $U_{f(x_0)}(\delta(f(x_0)))$, 
except that $s(g)\phi_{t_g}(D) \subset U_\lambda$ is replaced by $s(g)\phi_{t_f}(D) \subset U_{g(x_0)}(\delta(g(x_0)))$.) 

We note that 
${\cal F} \equiv cl\,{\cal E}^\#_{\delta/2}(X, M) \subset {\cal E}^\#_\delta(X, M)$, where the closure is taken in ${\cal E}^\#(X, M)$. 
In fact, each $f \in {\cal F}$ admits an open neighborhood ${\cal U}$ in ${\cal E}^\#(X, M)$ and a map $\Phi : {\cal U} \to {\cal H}(M)_0$ 
such that $\Phi(f) = id_M$ and $\Phi(g)f = g$ ($g \in {\cal U}$) (Proposition 2.5).
Each $g \in {\cal U} \cap {\cal E}^\#_{\delta/2}(X, M)$ ($\neq \emptyset$) admits an extension $\overline{g} \in {\cal E}^\#_{\delta/2}(D, M)$ 
and $\overline{f} = \Phi(g)^{-1} \overline{g}$ is an extension of $f$. 
If $g$ is sufficiently close to $f$, then $\Phi(g)$ is close to $id_M$ and $\overline{f} \in {\cal E}^\#_\delta(D, M)$. 
This means that $f \in {\cal E}^\#_\delta(X, M)$. 

Take a map $\mu : {\cal E}_{\cal U}^\#(X, M) \to [0, 1)$ with $\mu \geq \mu_1$ and $\mu|_{\cal F} \geq \mu_2|_{\cal F}$, 
and define the section $F$ by $F(f) = s(f) \phi_{\mu(f)}$. 

Finally we construct a f.p.\,deformation $\Psi_t$ ($t \in [0,1]$) of ${\cal E}_{\cal U}^\#(X, M)$ into ${\cal E}^\#_\delta(X, M)$ such that 
$\Psi_t({\cal E}^\#_{\delta/2}(X, M)) \subset {\cal E}^\#_\delta(X, M)$ ($0 \leq t \leq 1$). 
Using a cone structure of $D$ with the vertex $x_0$, we can find a s.d.r. $\psi_t$ ($t \in [0,1]$) of $D$ onto $x_0$ such that 
$\psi_t \in {\cal E}^+(D, D)$ ($0 \leq t < 1$) and $\psi_{t_1}(D) \supset \psi_{t_2}(D)$ ($0 \leq t_1 \leq t_2 \leq 1$).
If $\nu : {\cal E}_{\cal U}^\#(D, M) \to [0, 1)$ is sufficiently close to 1, 
then $h \psi_{\nu(h)} \in {\cal E}^\#_\delta(D, M)$ ($h \in {\cal E}_{\cal U}^\#(D, M)$) 
and $h \psi_t \in {\cal E}^\#_\delta(D, M)$ ($h \in {\cal E}^\#_\delta(D, M)$, $0 \leq t < 1$). 
The f.p. deformation $\Psi_t$ is defined by $\Psi_t(f) = F(f) \psi_{t \nu(F(f))}|_X$. 

It follows from (1) that the inclusion ${\cal E}^\#_{\delta/2}(X, M) \subset {\cal E}^\#_\delta(X, M)$ is a f.h.e.
Thus $\Psi_1 : {\cal E}_{\cal U}^\#(X, M) \to {\cal E}^\#_\delta(X, M)$ is a f.h. inverse of the inclusion 
${\cal E}^\#_\delta(X, M) \subset {\cal E}_{\cal U}^\#(D, M)$. 
This completes the proof of the case where $X$ is not an arc. 

When $X$ is an arc, 
there exists a s.d.r. $\psi_t$ of $X$ onto $x_0$ such that 
$\psi_t \in {\cal E}(X, X)$ ($0 \leq t < 1$) and $\phi_{t_1}(X) \supset \phi_{t_2}(X)$ ($0 \leq t_1 \leq t_2 \leq 1$).
Using $\psi_t$,  we can construct a f.p. deformation $\Psi_t$ ($t \in [0,1]$) of ${\cal E}_{\cal U}(X, M)$ into ${\cal E}_\delta(X, M)$ such that 
$\Psi_t({\cal E}^\#_\delta(X, M)) \subset {\cal E}^\#_\delta(X, M)$ ($0 \leq t \leq 1$). 

(2) The proof is similar to the arc case in (1) except that we use the s.d.r. $\psi_t(x) = (1 -t)x$ and the map $\nu_1(h) = \max\{ \nu(h), \nu(h\theta_{-1})\}$ 
instead of $\nu$ itself. 
\end{proof}

\begin{proof}[Proof of Proposition 6.7]
(1) In each case of (${\rm a}_\pm$) and (b), if ${\cal U_0}$ is the covering of $M$ by all open disks, then 
${\cal E}^\#(X, M) = {\cal E}_{\cal U_0}^\#(X, M)$ and ${\cal U}_\delta$ refines ${\cal U_0}$.  
Therefore, by Lemma 6.7, ${\cal E}_{\delta}^\#(X, M) \subset {\cal E}^\#(X, M)$ is a f.h.e. and 
${\cal E}_{\delta}(I, M) \subset {\cal E}(I, M)$ is a ${\Bbb Z}_2$-f.h.e. over $M$. 

(2) Any ${\cal U}$ admits a $\delta$ such that ${\cal U}_\delta$ refines ${\cal U_0}$. 
By (1) and Lemma 6.8 ${\cal E}_{\delta}^\#(X, M) \subset {\cal E}^\#(X, M)$ and ${\cal E}_{\delta}^\#(X, M) \subset {\cal E}_{\cal U}^\#(X, M)$ are f.h.e.'s. 
Thus the inclusion ${\cal E}_{\cal U}^\#(X, M) \subset {\cal E}^\#(X, M)$ is also a f.h.e.
Similarly, ${\cal E}_{\cal U}(I, M) \subset {\cal E}(I, M)$ is a ${\Bbb Z}_2$-f.h.e. over $M$. 
\end{proof}

%%%%%%%%%%%%%%%%%%%%%%%%%%%%%%%%%%%%%%%%%% Subsection 6.4.2 %%%%%%%%%%%%%%%%%%%%%%%%%%%%%%%%%%%%%%%%%%
\subsubsection{Proof of Theorem 6.1} \mbox{} 
\vskip 1mm
Finally, combining Propositions 6.6, 6.7 and Lemma 6.7, we can complete the proof of Propositions 6.1, 6.2 and Theorem 6.1

\begin{proof}[Proof of Proposition 6.1]
(1)(i) By Propositions 6.6\,(1)(i), 6.7\,(i)(a) and Lemma 6.7 we have the sequence of f.h.e.'s over $M$: 
\[ {\cal E}^\pm(X, M) \supset {\cal E}_\varepsilon^\pm(X, M) \cong {\cal E}_q^\pm(X, x; O_{TM}(\varepsilon), 0) \simeq S(TM). \]
(ii) The orientation  double cover $\tilde{M}$ has a canonical orientation. 
Let ${\cal U}$ be the open covering of $\tilde{M}$ consisiting of open disks $U$ on which $\pi : \tilde{M} \to M$ is injective. 
Each $f \in {\cal E}^\ast(X, M)$ admits a unique lift $\overline{f} \in {\cal E}^+_{\cal U}(X, \tilde{M})$ and 
this correspondence induces a f.p. homeomorphism ${\cal E}^\ast(X, M) \cong {\cal E}_{\cal U}^+(X, \tilde{M})$ over $M$. 
By (i) and Proposition 6.7\,(i)(a) we have the sequence of f.h.e.'s over $M$: 
\[ {\cal E}^\ast(X, M) \cong {\cal E}_{\cal U}^+(X, \tilde{M}) \subset {\cal E}^+(X, \tilde{M}) \simeq S(T\tilde{M}). \]
(2) By Propositions 6.6\,(2), 6.7\,(i) and Lemma 6.7 we have the sequence of f.h.e.'s over $M$: 
\[ {\cal E}(X, M) \supset {\cal E}_\varepsilon(X, M) \cong {\cal E}_q(X, x; O_{TM}(\varepsilon), 0) \cong S(TM). \]
\vskip -8.5mm
\end{proof}

\begin{proof}[Proof of Proposition 6.2] 
By Propositions 6.6\,(3), 6.7\,(ii) and Lemma 6.7 we have the sequence of ${\Bbb Z}_2$-f.h.e.'s over $M$:
\[ {\cal E}(I, M) \supset {\cal E}_\varepsilon(I, M) \cong {\cal E}_q(I, 0; O_{TM}(\varepsilon), 0) \simeq S(TM). \]
\vskip -8.5mm
\end{proof}

\begin{proof}[Proof of Theorem 6.1]
By Lemma 2.3 $X$ has a disk neighborhood $D$ in $M$. 
When $M$ is orientable, we orient $D$ and $M$ compatibly. 
From Proposition 6.1 it follows that \\
(i) if $M$ is orientable or $X$ is an arc, then ${\cal E}(X, M)_0 = {\cal E}^+(X, M) \simeq S(TM)$. \\ 
(ii) if $M$ is nonorientable and $X$ is not an arc, then ${\cal E}(X, M)_0 = {\cal E}^\ast(X, M) \simeq S(T\tilde{M})$. \\ 
\end{proof}

%%%%%%%%%%%%%%%%%%%%%%%%%%%%%%%%%%%%%%%%%% Reference %%%%%%%%%%%%%%%%%%%%%%%%%%%%%%%%%%%%%%%%%%

\end{document}